\documentclass[opre,nonblindrev]{informs3Custom} % current default for manuscript submission

\OneAndAHalfSpacedXI % current default line spacing
\usepackage{algorithm}
\usepackage{algorithmic}
\usepackage{multirow}
\usepackage{hhline}

\usepackage{endnotes}
\let\footnote=\endnote
\let\enotesize=\normalsize
\def\notesname{Endnotes}%
\def\makeenmark{$^{\theenmark}$}
\def\enoteformat{\rightskip0pt\leftskip0pt\parindent=1.75em
  \leavevmode\llap{\theenmark.\enskip}}

\usepackage{natbib}
 \bibpunct[, ]{(}{)}{,}{a}{}{,}%
 \def\bibfont{\small}%
 \def\bibsep{\smallskipamount}%
 \def\bibhang{24pt}%
\TheoremsNumberedThrough     % Preferred (Theorem 1, Lemma 1, Theorem 2)
\ECRepeatTheorems

\EquationsNumberedThrough    % Default: (1), (2), ...
\begin{document}
%%%%%%%%%%%%%%%%

%\setcounter{page}{1}
%\DOI{10.1287/xxxx.0000.0000}%

\RUNAUTHOR{Durante, Nascimento, and Powell}

\RUNTITLE{Backward ADP with HSMMs in Energy Storage Optimization}

% \TITLE{Bundling Information Goods of Decreasing Value}
\TITLE{Backward Approximate Dynamic Programming with Hidden Semi-Markov Stochastic Models in Energy Storage Optimization}

\ARTICLEAUTHORS{%
\AUTHOR{Joseph L. Durante}
\AFF{Department of Electrical Engineering, Princeton University, Princeton, NJ 08540, \EMAIL{jdurante@princeton.edu}} %, \URL{}}
\AUTHOR{Juliana Nascimento}
\AFF{Department of Operations Research, Princeton University, Princeton, NJ 08540, \EMAIL{jnascime@princeton.edu}}
\AUTHOR{Warren B. Powell}
\AFF{Department of Operations Research, Princeton University, Princeton, NJ 08540, \EMAIL{powell@princeton.edu}}
} % end of the block

\ABSTRACT{%
We consider an energy storage problem involving a wind farm with a forecasted power output, a stochastic load, an energy storage device, and a connection to the larger power grid with stochastic prices. Electricity prices and wind power forecast errors are modeled using a novel hidden semi-Markov model that accurately replicates not just the distribution of the errors, but also crossing times, capturing the amount of time each process stays above or below some benchmark such as the forecast. This is an important property of stochastic processes involved in storage problems. We show that we achieve more robust solutions using this model than when more common stochastic models are considered. The new model introduces some additional complexity to the problem as its information states are partially hidden, forming a partially observable Markov decision process. We derive a near-optimal time-dependent policy using backward approximate dynamic programming, which overcomes the computational hurdles of classical (exact) backward dynamic programming, with higher quality solutions than the more familiar forward approximate dynamic programming methods.
% Enter your abstract
}%

%\KEYWORDS{deterministic inventory theory; infinite linear programming duality;
%  existence of optimal policies; semi-Markov decision process; cyclic schedule}

\HISTORY{Version as of February 1, 2020}
\maketitle
%%%%%%%%%%%%%%%%%%%%%%%%%%%%%%%%%%%%%%%%%%%%%%%%%%%%%%%%%%%%%%%%%%%%%%

% Samples of sectioning (and labeling) in OPRE
% NOTE: (1) \section and \subsection do NOT end with a period
%       (2) \subsubsection and lower need end punctuation
%       (3) capitalization is as shown (title style).
%
%\section{Introduction.}\label{intro} %%1.
%\subsection{Duality and the Classical EOQ Problem.}\label{class-EOQ} %% 1.1.
%\subsection{Outline.}\label{outline1} %% 1.2.
%\subsubsection{Cyclic Schedules for the General Deterministic SMDP.}
%  \label{cyclic-schedules} %% 1.2.1
%\section{Problem Description.}\label{problemdescription} %% 2.

% Text of your paper here

\section{Introduction}
\label{Intro}

Renewable energy sources that exhibit high volatility and intermittency, such as wind or solar, are often paired with energy storage devices to improve the efficiency and reliability of the energy systems that incorporate them. To realize the full potential of the system, we must optimize the control policy to determine the best possible energy allocation decisions in an uncertain environment. Energy storage problems appear in many variations, characterized by the configuration of the system (which determines the controls) and the nature of the different types of uncertainties, such as prices, variability of wind and solar, and the behavior of the loads on the system.

In this paper we consider the problem of satisfying a time-varying load through a combination of energy from a storage device, a highly volatile wind power source, and the larger power grid with a highly volatile and heavy tailed locational marginal electricity price (LMP) at maximum profit. At every update of the LMP in the real time energy market (this occurs every 5 minutes in the Pennsylvania-New Jersey-Maryland Interconnection), energy allocation decisions between the nodes in this network must be made for the subsequent time interval. Specifically, we are able to: 1) send energy from the renewable source, storage device, or power grid to the load, 2) send energy from the renewable source to the storage device, or 3) buy or sell electricity using the connection between the storage device and the grid. We assume the role of a price taker and thus any energy flows between the system and the grid do not affect the LMP.

This is a common problem faced by owners of large industrial, residential, or commercial buildings equipped with a renewable source and an energy storage device. It is fairly standard to control these systems with simple rules, such as a buying energy when the electricity price is below some threshold and selling when it is above a different threshold. Such policies (known as policy function approximations) are often utilized as they are easy to implement, but it is quite difficult to adapt them to nonstationary environments, as is typical in energy applications. Deterministic lookahead approaches (also referred to as model predictive control or rolling horizon procedures) are also commonly used methods whose performance suffers in this highly stochastic environment. More sophisticated policies, such as those based on computing value functions for system states via dynamic programming or approximate dynamic programming, are necessary to maximize profits in this system.

However, even state-of-the-art dynamic programming methods are vulnerable to errors in the modeling of the underlying stochastic processes. A common pitfall when modeling energy storage problems is the use of standard time series models for wind power forecast errors, such as autoregressive (AR) models which are popular in the study of energy storage systems \citep[e.g.][]{lohndorf2010optimal, zhou2013managing}, that do not adequetly capture the crossing times of stochastic processes. This is a key characteristic of stochastic processes that arise in storage problems. Crossing times are contiguous blocks of time for which a stochastic process is above or below some reference series (in this case, a forecast). This paper models wind power generation with the univariate crossing state hidden semi-Markov model (HSMM) presented in \cite{durante2017} which is unique in its ability to capture both crossing time distributions and the distribution of errors from forecast. Additionally, the distribution of areas above and below the forecast (the surpluses and deficits of energy produced versus expected output) are accurately replicated by the crossing state model.

Characteristics such as crossing time, error, and area distributions are especially important to model in the context of an energy storage problem as they inform choices such as storage device capacity, type, and charge/discharge rate. Furthermore, properly modeling these behaviors influences the effectiveness of the control policy itself. Consider optimizing our energy system, whether through policy search, approximate dynamic programming (ADP), or another method, using a wind power model whose crossing times are too short compared to the true wind power production crossing time behavior (as is often the case when modeling wind power with an AR model). This will result in control policies that do not account for the possibility that wind will underperform expectations for extended periods of time. Thus, in practice, the resulting policy will not be robust as performance will likely suffer in these scenarios.

We extend the HSMM to model stochastic electricity prices in the real-time market by incorporating the temperature forecast as an explanatory variable. This produces very realistic sample paths of electricity prices, which can be particularly difficult to replicate given their heavy-tailed behavior and correlation with temperature.

The energy storage problem is formulated as a discrete time, finite horizon Markov decision process (MDP) in which a control decision must be made at each time step. In smaller, low dimensional problems, an optimal policy can be found using vanilla backward dynamic programming to compute value functions for each possible state. However, when a realistic system model is considered that incorporates more sophisticated stochastic models, the curses of dimensionality that arise in either the decision space, the state space, or the outcome space make performing a full backward pass computationally intractable. Specifically, when using HSMMs for the stochastics, the MDP becomes a partially observable MDP which introduces additional complexity into the problem. This motivates the need for an ADP method that will compute solutions efficiently without sacrificing much in terms of performance relative to the optimal solution (otherwise we lose the benefits of using a more sophisticated stochastic model).

This paper uses a special type of ADP known as {\it backward approximate dynamic programming}, as seen in \cite{senn2014reducing} and \cite{cheng2017low} for example, to overcome the curses of dimensionality. Unlike classical forward approximate dynamic programming, which estimates value functions while stepping forward in time (sometimes with a backward traversal), backward ADP performs a single backward pass, as done in standard backward dynamic programming, but then fits an approximate model based on a small sample of the states. In this setting, backward ADP has been found to produce near optimal solutions \citep[see][]{cheng2017low} which are much better than what was found using more familiar forward ADP methods \citep[see][]{jiang2014comparison}.

Furthermore, backward ADP and crossing state models are highly compatible as the models have a small (and fixed) number of compact post-decision information states for indexing the value function approximations. This allows for a reduction in both computation time and memory required to store value functions. By using backward ADP with the crossing state models, we can create a more realistic energy storage problem model and still find near-optimal control policies.

This paper makes the following contributions: 1) The exogenous processes involved in the co-located wind farm-energy storage device problem are modeled with a new hidden semi-Markov model that accurately replicates both crossing time and error distributions. 2) We explore the use of backward ADP, extending prior work by \cite{senn2014reducing} and \cite{cheng2017low} to higher dimensional problems, providing a robust complement to more classical forward ADP methods. 3) We extend the basic methodology of backward ADP to handle the hidden state variable of the hidden semi-Markov crossing state model, which required developing a more compact state representation, and the design of Bayesian updating logic. 4) We show that the backward ADP methodology produces higher quality results, more consistently, than forward ADP methods for energy storage problems, and further demonstrate that training on the hidden semi-Markov model produces more robust policies than standard stochastic models that have been used in the past.

The paper is organized as follows. Section \ref{Lit Rev} provides a brief literature review of energy storage problems, common stochastic models used in these problems, and algorithmic strategies related to backward ADP. Section \ref{HSMM} provides a thorough discussion regarding the modeling of the stochastic wind and price processes. Section \ref{Modeling} formally describes the energy storage problem by defining the five elements of the stochastic optimization problem. The backward ADP algorithms used in this paper are presented in Section \ref{backward ADP}. Numerical results comparing backward ADP to other policy types are reported in Section \ref{Numerical Results}, while results highlighting the impact of model selection on policy effectiveness and robustness are presented in Section \ref{Model Selection}. The paper is concluded in Section \ref{Conclusion}.

\section{Literature Review}
\label{Lit Rev}

%This paper first carefully modeling characteristics of the stochastics that are essential to the problem. The difficulties in controlling these renewable-based energy systems is that renewable power may deviate from its forecast than forecasted for extended periods of time and electricity prices may spike and then finds an optimization algorithm  involved in the energy storage problem first, and then develops 
Energy storage optimization is a widely researched topic with many problem variations. We review some of the problems and configurations that have been considered, organized by solution strategy. \cite{powell2016unified} describes four basic strategies for developing control policies for these systems: policy function approximations (PFAs), policies based on cost function approximations (CFAs), direct lookahead policies (DLAs), and policies based on value functions approximations (VFAs). PFAs map a state directly to a feasible action. CFAs maximize a parameterized approximation of a cost function subject to parametrically modified constraints. Both rely on policy search to optimize any parameters involved. DLAs maximize over both current and future actions based on an approximate model of the system; both deterministic and stochastic lookaheads fall in this class. VFAs maximize the one-step contribution of an action plus an approximation of the value of being in a future state. Extra attention is paid to VFA-based policies for the control of a single storage device in similar system configurations to the one considered in this paper. We then provide an overview of stochastic models used in dynamic programming approaches for energy storage control.

\subsection{Energy Storage Problem Variations and Solution Approaches}

%Before delving into the review, it is important to note that the effectiveness of each solution approach is highly dependent on the characteristics of the problem and how it is modeled. In general, the more complex and realistic the system model is, the more difficult it is to find optimal or near-optimal policies. For example, in the co-located wind farm-battery storage system considered in \cite{powell2016tutorial2}, a simple tuned buy-low, sell-high PFA will work best if the distribution of wind, electricity price, and load are assumed to be stationary. However, in reality there exists intertemporal correlations and daily patterns in the stochastic processes. If the system model incorporates these characteristics, the PFA will perform quite poorly, while policies from different classes will perform much better. Thus, multiple solution approaches are considered in this paper as we compare the policies resulting from backward ADP, the proposed solution method, to others from the same class (VFAs) and policies from other classes as well. 

Affine policies, a popular form of PFA, are linear functions that map states to actions. \cite{warrington2012robust}  and \cite{6275459} both use affine policies to control energy systems. \cite{han2016deep} train an artificial neural network (ANN), which is a form of nonlinear PFA, to optimize the flow of power between a wind farm, the grid, and a storage device to satisfy a time varying load. \cite{han2016deep} also demonstrate an ANN for a higher dimensional storage problem. 

A popular form of CFA is to use a deterministic lookahead (which can accommodate forecasts), with tunable parameters to handle uncertainty. This approach is used for robust power system control in \cite{SimaoReserves} where reserves are explicitly tuned to meet the variability of renewables. Similarly, in \cite{thalassinakis2004monte}, a Monte Carlo simulation method tunes the reserve level in addition to other power system settings. 

Deterministic lookaheads, also known as model predictive control \citep{camacho2013model}, are a popular policy for energy systems. In one example, \cite{denholm2009value} study the benefits of co-locating wind farms and and storage devices (using compressed air storage) using a deterministic lookahead which optimizes over a two-week rolling horizon. Model predictive control (MPC) is used by \cite{ma2012predictive} for optimizing heating and cooling systems in large buildings, while \cite{arnold2011model} use MPC to operate a storage hub with battery and hot water storage devices in the presence of uncertain renewable sources, electricity prices, and natural gas prices. \cite{4494596} optimize the bids of a wind farm and pumped storage facility using a stochastic lookahead, in the form of a two-stage stochastic program, to introduce uncertainty in the lookahead model. Two-stage stochastic programming is also utilized for much higher dimensional problems in the energy systems realm, such as robustly optimizing a large power grid in a day-ahead unit commitment problem \citep{wang2013two}.

%Stochastic lookaheads include both two-stage and multistage stochastic programming approaches. In \cite{4494596} two-stage stochastic programming is used to co-optimize bidding for a wind farm and pumped storage facility in the day-ahead market based on electricity price and wind power scenarios. The price distribution for each hour in the subsequent day is conditioned on the state of a hidden Markov model which allows for additional explanatory variables to influence the distribution. MULTISTAGE EXAMPLE

VFA-based policies approximate the impact of decisions now on the future through a statistical approximation of the value of being in a downstream state. \cite{schneider2015optimization} use ADP to optimize battery charging across a network of electric vehicle battery swap stations. Stochastic dual dynamic programming (SDDP) is a variant of ADP used for stochastic linear programs that has been applied to hydroelectric storage problems \citep{pereira1991multi, philpott2008convergence, shapiro2011analysis}. In this formulation, the stochastics involved are assumed to exhibit stage-wise (or intertemporal) independence. Other multistage stochastic programming algorithms such as the cut sharing algorithm from \cite{infanger1996cut} relax this assumption, allowing for interstage dependency.

An alternative approach for incorporating interstage dependencies is approximate dual dynamic programming (ADDP) which is used to optimize bidding and storage decisions in a hydro storage system in \cite{lohndorf2013optimizing}. \cite{lohndorf2015optimal} use ADDP to value natural gas storage and future trading, capturing interstage dependency by modeling forward prices as a discretized multivariate geometric Brownian motion (MGBM).

The MGBM model for futures prices is seen again in both \cite{lai2010approximate} and \cite{nadarajah2015relaxations} in which the management of commodity storage (such as natural gas) is considered. In \cite{lai2010approximate}, a novel forward ADP approach to valuating natural gas storage is used to benchmark commonly used heuristic valuation methods. In \cite{nadarajah2015relaxations}, VFAs for states in a commodity storage MDP are found using relaxations of approximate linear programs. While these gas storage problems are related to our battery storage problem, they are concerned with much longer time steps and optimization horizons (months to years) and experience different sources of uncertainty.

\subsection{Approximate Dynamic Programming for Energy Storage}

We now narrow our focus to the operation of a single energy storage device at finer timer scales. \cite{lohndorf2010optimal} combine approximate policy iteration with least squares policy evaluation (a form of forward ADP) to optimize the day-ahead bidding of a renewable supply-energy storage system participating in both the day-ahead and real-time electricity market, in which the price and renewable supply are modeled as correlated first order autoregressive processes.

\cite{zhou2013managing} utilize exact backward dynamic programming to find an optimal policy for controlling a co-located wind power-storage device system in the presence of stochastic wind and electricity price processes. Wind speed (transformed to wind power) is modeled as an AR(1) process with a seasonal component, while a carefully calibrated model with mean-reverting, seasonal, and jump components is used for prices. Results show that the exact dynamic programming approach results in considerable improvements over a rolling horizon procedure and other heuristic policies.

In \cite{cheng2017low}, a battery is co-optimized on different time scales for both energy arbitrage and frequency regulation using backward ADP. The highly correlated frequency regulation and LMP signal are modeled by forming ordered pairs of the prices and using a first order Markov chain to make hourly transitions between these ordered pairs. The LMP then evolves on a five minute time scale conditioned on its basepoint at the beginning of the hour.

Using the same wind farm-battery storage system configuration considered in this paper, \cite{jiang2014comparison} compare the effectiveness of several forward ADP methods based on approximate value iteration and approximate policy iteration in optimizing the system using different stochastic models of wind energy and electricity prices. Extensive testing shows that these classical ADP methods with several general purpose machine learning methods (linear models, tree regression, support vector regression, Gaussian process regression) work quite poorly in this energy storage problem, producing policies that achieve only 70-90 percent of optimality.

In contrast, backward ADP was found to produce very high quality results in realistic battery storage problems in \cite{cheng2017low}, with results ranging from 95 to 98 percent of optimality. Outside of the \cite{cheng2017low} paper, backward ADP has been relatively unexplored in the energy storage domain. \cite{senn2014reducing} also employs a backward ADP algorithm, although not in an energy storage application. Backward ADP is also closely related to numerical integration methods in dynamic programming that are prevalent in fields such as economics \citep[see][]{judd1998numerical,rust2008dynamic, cai2010stable}. The main difference is the method of randomly sampling states in the backward pass employed by the backward ADP algorithms, whereas numerical DPs often evaluate points in the state space on a grid before forming value function approximations. In addition, the numerical integration methods are typically limited to low-dimensional problems. Backward ADP is also highly related to approximate value iteration \citep[see, for example,][]{de2000existence, van2006performance, zanette2019limiting} and fitted value iteration \citep[see][]{munos2008finite}. It can be viewed as a special case of fitted value iteration in which the horizon is finite, one backwards iteration is performed, and the structure of the post-decision state is exploited to guide the sampling of the state space.

In this paper, we look to apply backward ADP to a higher dimensional problem. Backward ADP techniques from \cite{cheng2017low} are adapted to a setting where the stochastic wind and price processes are modeled with HSMMs. We explore both linear and lookup table approximation architectures. However, the backward ADP approach can be used with any parametric or nonparametric statistical learning methods \citep[see, for example,][]{hastie01statisticallearning}. The choice of approximation architecture is problem dependent, and should be chosen based on careful experimentation and benchmarking.

\subsection{Common Stochastic Models in Energy Storage Problems}

A common approximation for modeling wind (or the error in a wind forecast) is to use a first-order Markov process, whether an order-1 autoregressive process is used \citep[e.g.][]{lohndorf2010optimal, zhou2013managing} or a first-order Markov chain \citep[e.g.][]{jiang2014comparison, cheng2017low}. While the true wind process may be of a higher order, or depend on additional state-of-the-world variables, the information lost as a result of the simplification is often a necessary compromise to allow for the efficient computation of a solution when a dynamic programming approach is used. However, it may be the case that the simplified model is inaccurate with respect to the stochastic base model (or the actual behavior of the process) to the point that solution quality suffers in terms of expectation, robustness, or both.

Similar simplifying assumptions are made for electricity spot prices in dynamic programming approaches. For example, in \cite{sioshansi2014dynamic} perfect knowledge of future prices over the optimization horizon is assumed. A model with intertemporal independence is seen in \cite{xi2014stochastic} where the hourly spot price is determined by a lognormal distribution with parameters that depend on the hour of day. The first-order Markov process assumption is commonly used for prices as well, as is done with the mean reverting models used in \cite{tseng2002short} or \cite{mokrian2006stochastic}. This is also the case for the Markov model with jumps that is one of the models considered in \cite{jiang2014comparison} and the model with mean-reversion, jumps, and deterministic seasonality from \cite{zhou2013managing}.

%One common pitfall that causes some to overlook this fact is when policies are evaluated using sample paths that are generated with the simplified model to which value functions are fit instead of using the base model. On the other hand, this is done correctly in \cite{lohndorf2013optimizing} for example. Value functions are fit using a relaxed version of the problem, yet the policy is simulated forward with the stochastic base model. This is done in this paper as well -- value functions are found using a dimensionality-reduced problem in which the stochastics are simplified, but the testing is performed using the base model.

Most models commonly used in energy storage problems fail to replicate crossing time and error distributions, which are key characteristics of stochastic processes that can affect solution quality if not properly modeled. Meanwhile, the crossing state models described in \cite{durante2017} do so by construction. Furthermore, this can be accomplished with a relatively low-dimensional simplification of the information state variable. The result is a unique first-order Markov process that bridges the gap between simpler models which ignore important characteristics of the stochastics and complex models that require complicated information states (e.g. neural networks, models with many explanatory variables). Despite this modeling choice, finding exact solutions for the energy storage problems considered in this paper is still computationally expensive due to the dimensionality of the state, outcome, and decision space. This creates the need for backward ADP in this problem as it can compute high quality solutions efficiently.

\section{Modeling the Exogenous Processes}
\label{HSMM}

%The stochastic modeling of exogenous processes and, subsequently, the utilization of these models to develop more effective control policies are key contributions of this paper. Thus, modeling the stochastics will receive extra attention.

There are three exogenous processes that we have to model in this energy storage problem: the energy produced by wind, the electricity price, and the energy demand. Of these, the demand profile is far more predictable and exhibits far less variation from its expected value at each time $t$. In the interest of reducing the dimensionality of the outcome space and state space, we assume the load, $L_t$, follows a deterministic, time-dependent function formed from historical summertime demand profiles in the Princeton, NJ area. We consider three different demand profiles: one from a hot day, one from an average day, and one from a cool day. These all exhibit daily patterns such as troughs in the morning and peaks in the afternoon, but the total load is larger on hotter days. The series $\left\lbrace L_t \right\rbrace_{t=0}^T$ belongs in the initial state variable $S_0$ as it is a latent variable in our problem.

%At this point it is important to mention that throughout the paper we formulate our problem around decisions to send and receive quantities of energy (in kWh) to different nodes in the system over a fixed time interval. Load profiles, as well as wind forecasts and actual output, are reported as the instantaneous power (in kW) demand/production measurement at discrete points in time. To convert from instantaneous power to the energy demand/production over each time interval we simply multiply the instantaneous power by the length of the time step in hours. 

\subsection{The Crossing State Hidden Semi-Markov Model}
\label{WPM}

\begin{figure}
\centering
\includegraphics[width=5.5 in, height=2.3 in]{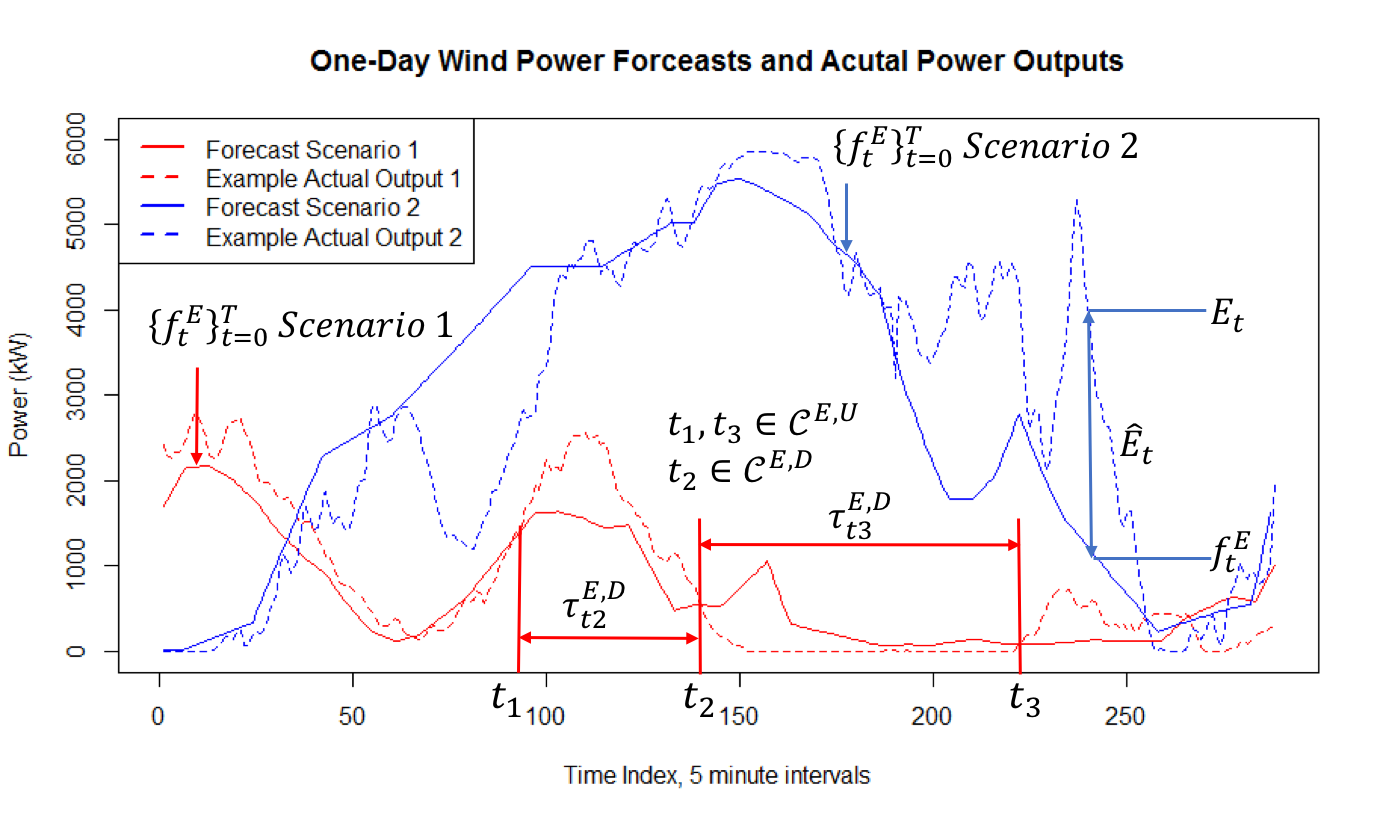}
\caption{Two wind forecast scenarios are shown (solid lines), along with corresponding actual wind power outputs (dotted). On the first forecast scenario, examples of points in time belonging to $\mathcal{C}^{E,U}$ and $\mathcal{C}^{E,D}$ and both a complete up- and down- crossing time are shown. On the second forecast path, a single forecast error, $\hat{E}_t=E_t-f_t^E$ is shown.}
\label{CrossingTimePic}
\end{figure}

\begin{figure}
\centering
\includegraphics[width=\columnwidth, height=2 in]{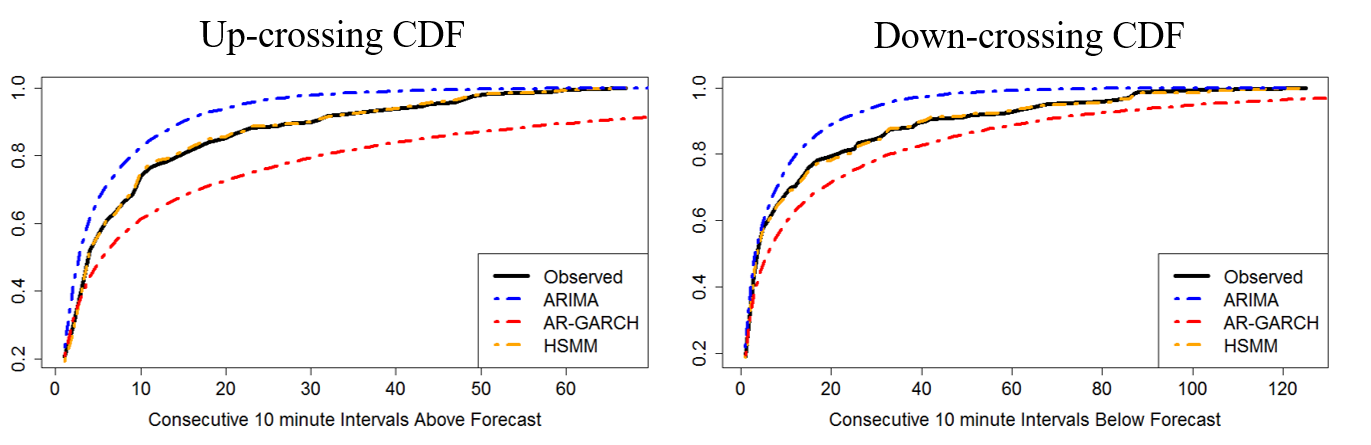}
\caption{Observed versus simulated up- (left) and down-crossing (right) time cumulative distribution functions for wind power forecast errors. The simulated distributions come from two common time series models, ARIMA and AR-GARCH, and the univariate crossing state HSMM presented in \cite{durante2017}. The crossing state model replicates crossing time distributions almost exactly (the CDF plots overlap the observed CDFs), while the time series model produce crossing times that are either consistently too long (AR-GARCH) or too short (ARIMA). This figure was adapted from \cite{durante2017}.}
\label{WindCrossingTime}
\end{figure}

Here we provide an overview of the univariate crossing state HSMM presented in \cite{durante2017}, including how to train the model. As mentioned in Section \ref{Intro}, we utilize this model for wind power production as it accurately captures the distribution of times for which wind power is above or below its forecast, otherwise known as the crossing times (examples are illustrated in Figure \ref{CrossingTimePic}). This behavior is poorly captured by standard time series models in comparison to the crossing state model, which replicates crossing time distributions for general univariate stochastic processes by construction \citep[see][]{durante2017}, as is shown in Figure \ref{WindCrossingTime}. In the model, we have conditional forecast errors distributions that are dependent on the partially hidden state variable, the crossing state. The ``semi-Markov" quality stems from the fact that crossing state transition probabilities are dependent on state duration.

Given a reference series, $\left\lbrace f^E_t \right\rbrace_{t=0}^T$, which in the case of wind power production is a power forecast, we define the error $\hat{E}_t=E_t-f^E_t$ where $E_t$ is the instantaneous wind power at time $t$. We assume a fixed forecast over the optimization horizon (it is a latent variable in the problem) and, therefore, $\left\lbrace f^E_t \right\rbrace_{t=0}^T \in S_0$. The wind power model is trained on sets of day-ahead wind power forecasts provided by an external vendor using a proprietary forecasting method and the resulting actual power output series. These are gathered from various wind farms in the Great Plains region; however, to form this univariate model, only data from the wind farm of interest is used.

First we define both up- and down-crossing times. Let the current elapsed time above forecast at time $t$, $\tau_{t}^{E,U}$,  be
\begin{align}
\tau_{t}^{E,U}=\ell \text{ if }
\begin{cases}
\hat{E}_{t-\ell} \leq 0\\
\hat{E}_{t+\ell'} > 0 \quad\forall \ell' \in \left\lbrace 0,1,...,\ell-1 \right\rbrace.
\end{cases}
\end{align} Similarly, the current elapsed time below forecast at time $t$, $\tau_{t}^{E,D}$, is defined as
\begin{align}
\tau_{t}^{E,D}=\ell \text{ if }
\begin{cases}
\hat{E}_{t-\ell} > 0 \\
\hat{E}_{t+\ell'} \leq 0 \quad\forall \ell' \in \left\lbrace 0,1,...,\ell-1 \right\rbrace.
\end{cases}
\end{align} Next, let the set of all indices such that forecast errors cross over from the negative to positive regime be $\mathcal{C}^{E,U}=\lbrace t|\hat{E}_{t-1}\leq 0 \wedge \hat{E}_t >0\rbrace$. Likewise, let the set of all indices such that errors cross over from the positive to negative regime be $\mathcal{C}^{E,D}=\lbrace t|\hat{E}_{t-1}\geq 0 \wedge \hat{E}_t <0\rbrace$. The sets of all up- and down-crossings are then $\mathcal{T}^{E,U}=\lbrace \tau_{t}^{E,U}|t+1 \in \mathcal{C}^{E,D}\rbrace$ and $\mathcal{T}^{E,D}=\lbrace \tau_{t}^{E,D}|t+1 \in \mathcal{C}^{E,U}\rbrace$ respectively as a crossing time is simply a completed elapsed time above or below forecast. Examples of points in time belonging to $\mathcal{C}^{E,U}$ and $\mathcal{C}^{E,D}$, as well as up- and down- crossing times are shown in Figure \ref{CrossingTimePic}.

For both up- and down-crossing times, there exists cumulative distribution functions $F^{E,U}$ and $F^{E,D}$ respectively. Up-crossing time distributions are quantized by partitioning into $m^E$ bins, splitting at the $q_i=\frac{i}{m^E}$ quantile points for $i=0,1,...,m^E-1$. An up-crossing time, $\tau_t^{E,U} \in \mathcal{T}^{E,U}$, belongs to crossing time duration bin $D_t^{E}=d$ if $q_d \leq F^{E,U}(\tau_t^{E,U}) < q_{d+1}$. Down-crossing time distributions are similarly quantized.

Our crossing state variable $I_t^{C,E}\equiv (C_t^{E},D_t^{E})$ is defined as the pair of variables describing whether or not the error is above the forecast, $C_t^{E}=\mathbf{1}_{\left\lbrace \hat{E}_t>0\right\rbrace}$, and to which crossing time duration bin, $D_t^{E}$, the completed crossing time will belong. Note that this means that during online optimization, the state $C^E_t$ is observable (we know if we are above or below the forecast), but the duration bin $D_t^E$ is not until the sign variable $C^E_t$ switches. However, when building the crossing state-dependent error distributions from training data for the model we can observe the duration bin at each point in time by peeking into the future to find the complete crossing time. Letting $\mathcal{I}^{C,E}$ be the set of all possible crossing states for the process, there exists a distribution from training data of crossing times $F^{\tau,E}_{i}$ for each possible crossing state $i \in \mathcal{I}^{C,E}$. These distributions serve as the sojourn time distributions for the crossing states.

Transitions between crossing states are made using a transition matrix $\mathbb{P}(i'|i)$ for each pair of crossing states $(i',i) \in \mathcal{I}^{C,E} \times \mathcal{I}^{C,E}$ in which self-transitions are not allowed ($\mathbb{P}(i|i)=0$ $\forall i \in \mathcal{I}^{C,E}$). This matrix is computed from training data by considering only pairs of points in time $(t,t+1)$ such that $t+1 \in (\mathcal{C}^{E,U} \cup \mathcal{C}^{E,D})$ (points where the crossing state makes a transition). For all of these pairs, letting $n(I_{t+1}^{C,E}=i'|I_t^{C,E}=i)$ be the count of the transitions from state $i$ to state $i'$ occurring for each pair of crossing states $(i,i')$ and $n(I_t^{C,E}=i)$ be the number of times $I_t^{C,E}=i$ for each crossing state $i$, the empirical transition probability from crossing state $i$ to $i'$ is
\begin{align}
\label{tprob}
\mathbb{P}(i'|i)=\frac{n(I_{t+1}^{C,E}=i'|I_t^{C,E}=i)}{n(I_t^{C,E}=i)}.
\end{align} The duration-dependent crossing state transition probability is then a function of the elapsed time above or below forecast given by
\begin{align}
\mathbf{P}(I_{t+1}^{C,E}=i'|I_{t}^{C,E}=i,\tau^E_t)=
\begin{cases}
1-F_{i}^{\tau,E}(\tau_t^E) & \text{if } i'=i\\
F_{i}^{\tau,E}(\tau_t^E)\mathbb{P}(i'|i) & \text{for } i'\neq i.
\end{cases}
\end{align}

\begin{figure}
\centering
\includegraphics[width=\columnwidth, height=2 in]{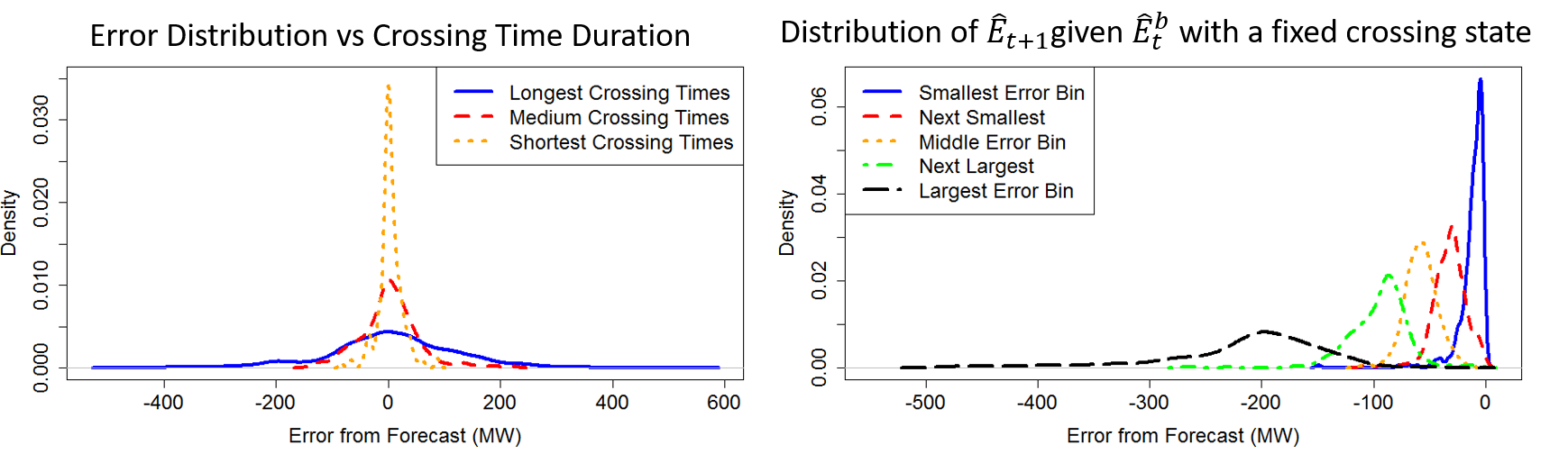}
\caption{Left: Error distributions conditioned on $D_t^E$ with $m^E=3$. Both positive and negative errors for equal values of $D_t^E$ are combined to form the distributions. The variance of the errors tends to increase with run length. Right: Example of conditional distributions for $\hat{E}_{t+1}$ given a fixed crossing state, $I_t^{C,E}=(0,2)$, but varying which error bin, $\hat{E}_t^b$, that $\hat{E}_t$ belongs to. The magnitude of the next error is largely dependent on the magnitude of the current error. This figure is a slightly modified version of a figure from \cite{durante2017}.}
\label{RLErrors}
\end{figure}

Next we describe conditioning the error generation on the crossing state. From training data, there exists empirical conditional error CDFs $F^{\hat{E}}_i$ and corresponding error density functions $\mathbf{P}(\hat{E}_{t+1}|i)$ for $i \in \mathcal{I}^{C,E}$. Error distributions are not identical across crossing states; in fact they are likely to be quite different, such as the case where the error distribution is asymmetric. Furthermore, error distributions are likely to vary across duration bins as well. This behavior is seen in the left plot of Figure \ref{RLErrors}, which shows error densities for each type of duration bin with $m^E=3$. Thus, to better capture the behavior of the error process, the error generation process is conditioned on the crossing state $I_t^{C,E}$.

In addition to errors being crossing state-dependent, they are dependent on error history as well. A first order Markov chain is used to model this behavior. Similar to how the crossing time distributions are quantized, each conditional error distribution $F^{\hat{E}}_i$ is partitioned into $n^E$ bins, splitting at the $q_j=\frac{j}{n^E}$ quantile points for $b=0,1,...,n^E-1$. The error $\hat{E}_t$ belongs to bin $\hat{E}^b_t$ if $q_b \leq F^{\hat{E}}_i(\hat{E}_t) < q_{b+1}$. Then, given $\hat{E}_t \in \hat{E}^b_t$, we form conditional empirical distributions for the error at time $t+1$ giving $\mathbf{P}(\hat{E}_{t+1}|I_t^{C,E}, \hat{E}_t^b)$. The dependence of $\hat{E}_{t+1}$ on $E_t^b$, the aggregated state of the current error, is illustrated in the right plot of Figure \ref{RLErrors} in which conditional distributions for $\hat{E}_{t+1}$ are plotted for a fixed crossing state, but varying error states $\hat{E}_t^b$.

It is important to realize that the same error $\hat{E}_t$ can fall in different error bins for different crossing states. For example, the error $\hat{E}_t=+2000$ kW may be in bin $\hat{E}_t^b=4$ for the $I_t^{C,E}=(1,0)$ crossing state (short up-crossings), but for the $I_t^{C,E}=(1,2)$ state (longer up-crossings), it may belong to bin $\hat{E}_t^b=2$. To avoid additional notation, \textit{the variable $\hat{E}_t^b$ will always be paired with a crossing state and refers to the bin that the error $\hat{E}_t$ belongs to for the corresponding crossing state}.

For each crossing state $i \in \mathcal{I}^{C,E}$, there also exists an error density $\mathbf{P}(\hat{E}_{t+1}|i, t+1 \in \mathcal{C}^{E,U} \cup \mathcal{C}^{E,D})$. This is the distribution of the initial error given the process has just transitioned to the new crossing state $i$.

The information state variable for this process, denoted $I_t^E$, contains the following variables at each time $t$: $I_t^E \equiv \left( C_t^E,D_t^E,\tau_t^E,\hat{E}^b_t \right)\equiv \left(I_t^{C,E},\tau_t^E,\hat{E}^b_t \right)$. If known, these variables fully determine the distribution of the exogenous information $\hat{E}_{t+1}$, given by
\begin{align}
\label{Ehat}
\begin{split}
\mathbf{P}(\hat{E}_{t+1}|I_t^{C,E}=i,\tau_t^E,\hat{E}^b_t)=&(1-F_{i}^{\tau,E}(\tau_t^E))\mathbf{P}(\hat{E}_{t+1}|i,\hat{E}^b_t)
+\\
&F_{i}^{\tau,E}(\tau_t^E)\sum\limits_{i' \neq i}\mathbb{P}(i'|i)\mathbf{P}(\hat{E}_{t+1}|i', t+1 \in \mathcal{C}^{E,U} \cup \mathcal{C}^{E,D}).
\end{split}
\end{align}

Finally, we need to model the price charged by the grid in the real-time market, $P_t$, a quantity known as the "Locational Marginal Price" (or LMP). Electricity prices exhibit extreme spikes, often in excess of 20 or 30 times the base price, as well as a dependence on temperature, and a daily cyclic pattern. We extend the above crossing state model to capture these characteristics by including a seasonal (daily) periodic component and conditioning price distributions on a temperature forecast. To differentiate the price model from the wind model, $P$ takes the place of every $E$ in the notation.

A thorough discussion of the price model including implementation details is given in the Appendix, but the model is summarized here. An LMP forecast $f_t^P=\mathbb{E}[P_t] \text{ } \forall t$ is assumed to be deterministic over the optimization horizon and thus $\left\lbrace f^P_t \right\rbrace_{t=0}^T \in S_0$. Additionally, we isolate two variables at each point in time $t$ from an input temperature forecast: $y^s_{t}$, a quantized seasonality component of the temperature forecast, and $ y^{tr}_{t}$, a quantized trend component of the temperature forecast. These also belong in the initial state variable. We then condition the error-from-reference densities on these additional variables, forming the distributions $\mathbf{P}(\hat{P}_{t+1}|I_t^{P}, y^s_{t+1}, y^{tr}_{t+1})$ for each information state $I_t^P$.

\subsection{Compact Exogenous Information States}
\label{CompactInfo}

Letting $\tau_i^{E,max}$ be the largest crossing time for crossing state $i\in \mathcal{I}^{C,E}$, we see that there are $\sum\limits_{i\in \mathcal{I}^{C,E}} n^E \tau_i^{E,max}$ possible information states at each time $t$. This number can be quite large, especially if crossing times tend to span many time periods. Fitting value functions to system states with backward ADP will likely be too computationally expensive without a more compact information state representation. For this reason, we introduce a modified compact process information state $\tilde{I}_t^E \equiv\left(C_t^E,D_t^E,\hat{E}^b_t \right)$, which can take on $2 \times m^E \times n^E$ states. Similarly, the compact information state representation for the electricity price process is given by $\tilde{I}_t^P \equiv\left(C_t^P,D_t^P,\hat{P}^b_t \right)$.

Note that the error distributions $\mathbf{P}(\hat{E}_{t+1}|i,\hat{E}^b_t)$ and $\mathbf{P}(\hat{E}_{t+1}|i,t+1 \in \mathcal{C}^{E,U} \cup \mathcal{C}^{E,D})$ for all $i \in \mathcal{I}^{C,E}$ are unaffected by this change. However, the transition between crossing states must now be modeled with a Markov approximation of the semi-Markov model as no elapsed time above or below forecast count is maintained. Transition probabilities are now given by a time-invariant modified crossing state transition matrix $\tilde{\mathbb{P}}(I_{t+1}^{C,E}|I_t^{C,E})$ which allows for self-transitions. This is estimated from training data using Equation \ref{tprob}; however all time periods $t$ are considered, not only pairs of points where errors switch signs. 

%Consequently, for each crossing state $i$, a Bernoulli random variable with probability of success given by $\tilde{\mathbb{P}}(i|i)$ determines whether or not the process remains in state $i$ between time periods. As a result, the crossing time in each crossing state is a realization of a $Geo(\tilde{\mathbb{P}}(i|i))$ random variable.  Also note that for $i'\neq i$, the semi-Markov transition probability can be related to the Markov model transition probability by $\mathbb{P}(i'|i)=\tilde{\mathbb{P}}(i'|i)/(1-\tilde{\mathbb{P}}(i|i))$.

Consequently, the distribution of $\hat{E}_{t+1}$, given only the compact information state, is
\begin{align}
\mathbf{P}(\hat{E}_{t+1}|I_t^{C,E}=i,\hat{E}^b_t)=\tilde{\mathbb{P}}(i|i)\mathbf{P}(\hat{E}_{t+1}|i,\hat{E}^b_t)
+\sum\limits_{i' \neq i}\tilde{\mathbb{P}}(i'|i)\mathbf{P}(\hat{E}_{t+1}|i', t+1 \in \mathcal{C}^{E,U} \cup \mathcal{C}^{E,D}).
\end{align} We will be fitting value functions using these compact information states with backward ADP and the Markov approximation for the process, while using the full information state and the semi-Markov model in forward passes. 

\subsection{The Belief State and its Bayesian Update}
\label{Bayesian Update}

While observing the above stochastic processes in the forward pass, the system operator will know the magnitude and sign of the current error from the reference series for each process as well as both $\tau_t^E$ and $\tau_t^P$, the elapsed times above or below forecast. However, the crossing time duration bins $D_t^E$ and $D_t^P$ are unknown at time $t$. Consequently, $\hat{E}_t^b$ and $\hat{P}_t^b$ are unknown as well. Thus, the information states are partially unobservable.

As a VFA-based policy takes the action that maximizes the one-step contribution plus the expected value of the downstream state, we must know the probability of reaching each downstream state when taking an action. This requires \textit{belief states}, denoted $B_t^{E}$ and $B_t^{P}$ for wind and price respectively, giving the operator's distribution of belief about the unknown variables. The remainder of this section discusses the belief state and Bayesian update for the wind process, but note that there are analogous time $t$ beliefs and Bayesian update functions for the price process.

Let $B_t^E \equiv \left(\left\lbrace \mathbf{P}(I_t^{C,E}=i)\right\rbrace_{i\in\mathcal{I}^{C,E}},\tau_t^E, \hat{E}_t \right)$. Given $B_t^E$, we can derive our belief about the distribution of the error at time $t+1$, $\mathbf{P}(\hat{E}_{t+1}|B_t^E)$. We are able to determine the sign of the error, $C_t^E$, based on $\hat{E}_t$. Then, for each possible value of $D_t^E$, and corresponding crossing state $I_t^{C,E}=(C_t^E,D_t^E)$, $\hat{E}_t$ can belong to only one error bin $\hat{E}_t^b$. We then have
\begin{align}
\label{Etplus1givenKt}
\mathbf{P}(\hat{E}_{t+1}|B_t^E)=\sum\limits_{i\in \mathcal{I}^{C,E}} \mathbf{P}(I_t^{C,E}=i)\mathbf{P}(\hat{E}_{t+1}|i,\tau_t^E,\hat{E}^b_t),
\end{align}
where $\mathbf{P}(\hat{E}_{t+1}|i,\tau_t^E,\hat{E}^b_t)$ is given by Equation \ref{Ehat}.

Subsequently, following the observation of $\hat{E}_{t+1}$, a Bayesian update is performed on the belief state at each time step according to the update function $B_{t+1}^E=U^{E}(B_t^E,\hat{E}_{t+1})$ defined by the following two cases:

\begin{hitemize}
\item Case 1: $sign(\hat{E}_{t+1})=sign(\hat{E}_{t})$. In this case $\tau_{t+1}^E=\tau_{t}^E+1$. This is then used to compute the likelihood that the future completed crossing time belongs to bin $D_t^E$ for each crossing state $i=(C_t^E,D_t^E)$ given the elapsed time above or below forecast $\tau^E_{t+1}$. This likelihood is given by $1-F^{\tau,E}_{i}(\tau^E_{t+1})$. Furthermore, the likelihood of observing error $\hat{E}_{t+1}$ in crossing state $i$ given a crossing state transition has not occurred is $\mathbf{P}(\hat{E}_{t+1}|i,\hat{E}^b_t)$. Thus with prior beliefs, $\mathbf{P}(I_{t}^{C,E}=i)$ for $i\in\mathcal{I}^{C,E}$, we compute the posterior beliefs
\begin{equation}
\mathbf{P}(I_{t+1}^{C,E}=i)=\frac{1}{p^{norm}}\left(\mathbf{P}(I_{t}^{C,E}=i)(1-F^{\tau,E}_{i}(\tau^E_{t+1}))\mathbf{P}(\hat{E}_{t+1}|i,\hat{E}^b_t)\right),
\end{equation} where $p^{norm}=\sum\limits_{i'\in\mathcal{I}^{C,E}}\mathbf{P}(I_{t}^{C,E}=i')(1-F^{\tau,E}_{i'}(\tau^E_{t+1}))\mathbf{P}(\hat{E}_{t+1}|i',\hat{E}^b_t)$.

\item Case 2: $sign(\hat{E}_{t+1})\neq sign(\hat{E}_t)$. In this case $\tau_{t+1}^E=0$ and we are able to determine the crossing state at time $t$ based on the sign of $\hat{E}_t$ and the completed crossing time $\tau^E_t$; let this be state $i^*$. Additionally, we know that a crossing state transition has taken place. The likelihood of observing error $\hat{E}_{t+1}$ in crossing state $i$ given a crossing state transition has just occurred is $\mathbf{P}(\hat{E}_{t+1}|i, t+1 \in \mathcal{C}^{E,U} \cup \mathcal{C}^{E,D})$. Thus, for $i\in\mathcal{I}^{C,E}$, posterior beliefs are given by
\begin{equation*}
\mathbf{P}(I_{t+1}^{C,E}=i)=\frac{1}{p^{norm}}\left(\mathbb{P}(i|i^*)\mathbf{P}(\hat{E}_{t+1}|i, t+1 \in \mathcal{C}^{E,U} \cup \mathcal{C}^{E,D})\right),
\end{equation*} where $\mathbb{P}(i|i^*)$ is defined in Equation \eqref{tprob} and $p^{norm}=\sum\limits_{i'\in\mathcal{I}^{C,E}}\mathbb{P}(i'|i^*)\mathbf{P}(\hat{E}_{t+1}|i',t+1 \in \mathcal{C}^{E,U} \cup \mathcal{C}^{E,D})$.
\end{hitemize} Given these recursive updating formulas for the belief state, we only need to initialize our beliefs at $t=0$. Given $\hat{E}_0$, we use a discrete uniform distribution for the initial beliefs: $\mathbf{P}(I_0^{C,E}=i)=1/(m^E)$ for $i \in \mathcal{I}^{C,E}$ such that $C_t^E=sign(\hat{E}_0)$. Setting $\tau^E_0=0$, this forms $B_0^E$.

\section{Mathematical Model of the Energy Storage Problem}
\label{Modeling}

\begin{figure}
\centering
\includegraphics[width=4.5 in, height=1.6 in]{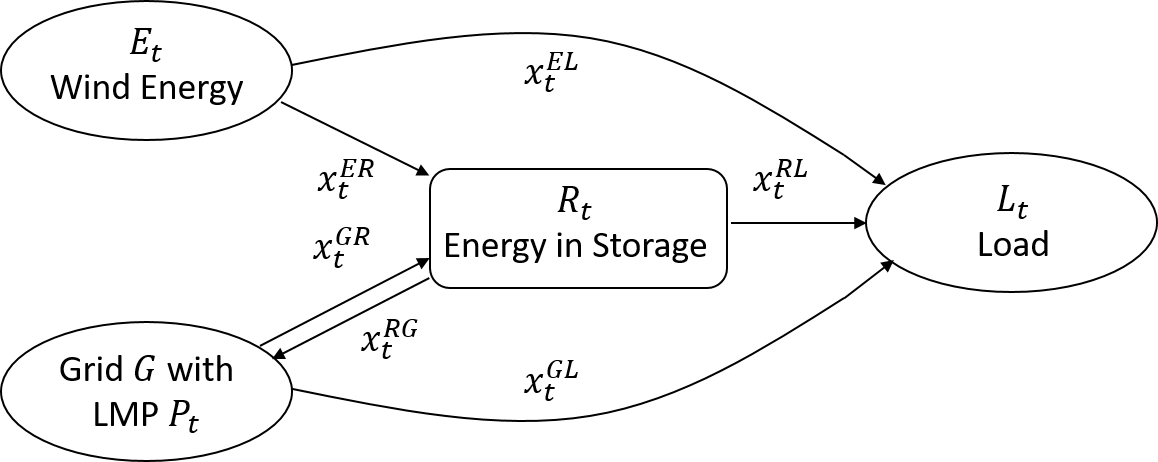}
\caption{An energy storage problem with four energy nodes: wind, the larger power grid, storage, and demand/load; three exogenous variables: wind energy production, the price of electricity, and the demand; and six decision variables representing possible energy allocations at each time step.}
\label{EnergyStorageProblem}
\end{figure}

We consider an energy storage problem, similar to the model presented in \cite{powell2016tutorial2}, in which a stochastic wind energy supply, a battery, and the power grid, with an associated stochastic electricity price, are used in combination to satisfy a time-varying power demand. The objective is to control the system, whose configuration is illustrated in Figure \ref{EnergyStorageProblem}, at minimum cost (or maximum profit). There are four nodes, six decision variables, and three exogenous processes in this energy storage system configuration.

This section provides a complete model of the stochastic optimization problem by defining the state variable, the decision variable, exogenous information, the transition function, and the objective function. The modeling style and notation are adopted from \cite{powell2016unified} and \cite{powell2016tutorial2}.

\subsection{The State Variable}

First, we define the dynamic state variable $S_t$ and the initial state variable $S_0$. Following the modeling convention from \cite{powell2016unified}, the initial state contains all data pertinent to the problem, including constants and deterministic variables. Conversely, $S_t$ contains only variables which may change over time.

We are treating forecasts as static in this problem, thus they belong in $S_0$. In addition to latent variables, $S_0$ contains initial beliefs about unknown parameters and initial values of wind, price, and storage level. Thus, our initial state is given by
\begin{align}
S_0=\left(\left\lbrace f^E_t, f^P_t, L_t, y^s_t, y^{tr}_t \right\rbrace_{t=0}^T, R^{max}, \eta,\rho^{ch},\rho^{dch}, R_0, E_0, P_0, B_0^E, B_0^P \right),
\end{align} where $R^{max}$ is the maximum capacity of the battery, $\eta \in [0,1]$ is the battery round trip efficiency, and $\rho^{ch}$ and $\rho^{dch}$ are maximum battery charge and discharge rates respectively.

For our VFA-based policies, we must track our time $t$ beliefs in the belief states $B_t^{E}$ and $B_t^{P}$. The contribution function, $C(S_t,x_t)=P_t(L_t-x_t^{GR}-x_t^{GL}+\eta x_t^{RG})$, requires that the current electricity price be known, thus $P_t=\hat{P}_t+f_t^P$ is incorporated as well. Finally, $E_t=\hat{E}_t+f_t^E$ and $R_t$ are included to determine our constraints on the decision vector $x_t$. Thus, our dynamic state variable $S_t$ for $t>0$ is given by
\begin{align}
S_t=(R_t, E_t, P_t, B_t^E, B_t^P).
\end{align}

Note that in this problem, the post-decision state variable $S_t^x$, which carries only the information necessary to transition to $S_{t+1}$ \textit{after} a decision has been made \citep{powell2011approximate}, is given by
\begin{align}
S_t^x=(R_t^x, B_t^{E}, B_t^{P}).
\end{align} $R_t^x$ is the energy level of the battery following the decision to use or store energy at time $t$.

\subsection{The Decision Variable}
The decision variable at each time $t$ is
\begin{align}
x_t=(x_t^{GL},x_t^{GR},x_t^{RG},x_t^{EL},x_t^{ER},x_t^{RL}),
\end{align} where $x_t^{AB}$ indicates energy sent from node $A$ to node $B$. This is subject to the constraints
\begin{align}
x_t^{EL}+x_t^{ER} &\leq E_t,\label{D1}\\
x_t^{GL}+x_t^{EL}+\eta x_t^{RL} &= L_t,\label{D2}\\
x_t^{RG}+x_t^{RL} &\leq \min(R_t,\rho^{dch}),\label{D3}\\
x_t^{ER}+x_t^{GR} &\leq \min(\rho^{ch},R^{max}-R_t),\label{D4}\\
x_t^{GL},x_t^{GR},x_t^{RG},x_t^{EL},x_t^{ER},x_t^{RL} &\geq 0\label{D5}.
\end{align} All vector decisions $x_t$ that satisfy the above constraints form the set of feasible decisions at time $t$, $\mathcal{X}_t(S_t)$, or $\mathcal{X}_t$ where the dependence on $S_t$ is implied. However, for the dynamic programming approaches in this paper, we assume that $x_t^{EL}=\min(E_t,L_t)$, and $x_t^{ER}=\min\left(\rho^{ch},R^{max}-R_t,E_t-x_t^{EL}\right)$ to reduce the size of the feasible region. Additionally, we assume any unused wind energy $E_t-(x_t^{EL}+x_t^{ER})$ is dissipated.

Constraint \eqref{D1} ensures that the amount of renewable energy used at time $t$ does not exceed the amount produced. Constraint \eqref{D2} requires that the load be met at each time $t$. Limits on the amount of energy that can be drawn from and sent to the battery at each time step are imposed in constraints \eqref{D3} and \eqref{D4} respectively. Finally, a non-negativity constraint is imposed on each element of the decision vector in constraint \eqref{D5}. Note that there is no constraint on the amount of energy that can be purchased from the grid as we assume it is an infinite source of power.

\subsection{Exogenous Information}
\label{Exo Info}

Section \ref{HSMM} was devoted to the stochastic modeling of the exogenous wind and electricity price processes. The exogenous information arriving between $t$ and $t+1$ is given by $W_{t+1}=(\hat{E}_{t+1}, \hat{P}_{t+1})$ where our time $t$ belief about the distribution of $\hat{E}_{t+1}$ is given by $\mathbf{P}(\hat{E}_{t+1}|B_t^E)$ and that of $\hat{P}_{t+1}$ is given by $\mathbf{P}(\hat{P}_{t+1}|B_t^P)$. It is important to note that in our model we assume we are a price taker and any decisions made do not affect the electricity price process. 

\subsection{The Transition Function}
\label{Transition Function}

The system transition function, $S_{t+1}=S^M(S_t,x_t,W_{t+1})$, is defined by
\begin{align}
& R_{t+1}=R_t+\eta(x_t^{GR}+x_t^{ER})-x_t^{RL}-x_t^{RG}, \label{Trans 1}\\
& E_{t+1}=f_{t+1}^E+\hat{E}_{t+1}, \label{Trans 2}\\
& P_{t+1}=f_{t+1}^P+\hat{P}_{t+1}, \label{Trans 3}\\
& B_{t+1}^E=U^{E}(B_t^E,\hat{E}_{t+1}), \label{Trans 7}\\
& B_{t+1}^P=U^{P}(B_t^P,\hat{P}_{t+1}), \label{Trans 11}
\end{align} where $U^{E}(B_t^E,\hat{E}_{t+1})$ and $U^{P}(B_t^P,\hat{P}_{t+1})$ are the Bayesian updating functions defined in Section \ref{Bayesian Update}. Note this can be broken into two functions: the pre- to post- decision state transition function, $S_t^x=S^{M,x}(S_t,x_t)$, and the transition function from post-decision state to the next pre-decision state given the arrival of exogenous information $W_{t+1}$, $S_{t+1}=S^{M,W}(S_t^x,W_{t+1})$. $S_t^x=S^{M,x}(S_t,x_t)$ alters the resource energy level based on the decision $x_t$ as in Equation \eqref{Trans 1}: $R^x_{t}=R_t+\eta(x_t^{GR}+x_t^{ER})-x_t^{RL}-x_t^{RG}$. Additionally, $P_t$ and $E_t$ are dropped from pre- to post- decision state, and the remaining variables remain unaltered. $S_{t+1}=S^{M,W}(S_t^x,W_{t+1})$ is then given by Equations \eqref{Trans 2}-\eqref{Trans 11} along with $R_{t+1}=R_t^x$.

\subsection{The Objective Function}
\label{Objective Function}

As we aim to operate the system at maximum profit in the real-time electricity market, each time step represents a five-minute interval in our program, and our finite horizon control problem has the objective function
\begin{align}
\operatornamewithlimits{max}\limits_{\pi \in \Pi}\mathbb{E}^\pi \left[\sum\limits_{t=0}^T C(S_t,X^{\pi}_t(S_t))|S_0\right],
\end{align}
where the contribution function is given by $C(S_t,x_t)=P_t(L_t-x_t^{GR}-x_t^{GL}+\eta x_t^{RG})$ and $S_{t+1}=S^M(S_t,X^{\pi}_t(S_t),W_{t+1})$ is determined by the system transition function. We are maximizing over the set of all possible policies $\pi \in \Pi$. In the contribution function we profit from satisfying the load or selling energy back to the grid, but must pay for any energy that originates from the grid.

\section{Backward Approximate Dynamic Programming}
\label{backward ADP}

Assuming a terminal reward of $V_T^*(S_T)$ for each terminal state $S_T$, if computationally tractable, an optimal policy can be found using vanilla backward dynamic programming to find value functions, $V^*_{t}(S_t)$, for each possible system state. Value functions are given by Bellman's equation for finite horizon problems,
\begin{align}
\label{ExactVal}
V_t^*(S_t)=\max\limits_{x_t \in \mathcal{X}_t} \left(C(S_t,x_t)+\mathbb{E}\left[V^*_{t+1}(S_{t+1})|S_t,x_t\right]\right),
\end{align}
and, once these are found, the optimal policy, 
\begin{align}
\label{ExactPol}
X_t^{*}(S_t)=\argmax\limits_{x_t \in \mathcal{X}_t} \left(C(S_t,x_t)+\mathbb{E}\left[V^*_{t+1}(S_{t+1})|S_t,x_t\right]\right),
\end{align} maximizes the one-step contribution plus the expected value of the downstream state. In cases where performing a complete backward pass is either impossible or impractical (as is the case in our problem), we can instead rely on approximations of these value functions and use a VFA-based policy,
\begin{align}
\label{ApproxPol}
X_t^{\pi}(S_t)=\argmax\limits_{x_t \in \mathcal{X}_t} \left(C(S_t,x_t)+\mathbb{E}\left[\bar{V}_{t+1}(S_{t+1})|S_t,x_t\right]\right),
\end{align} where $\bar{V}_{t+1}(S_{t+1})$ is some approximation of the value of the downstream states.

To remove the expectation from the policy, we can fit value functions instead to the post-decision state variable $S_t^x$. This is possible for our problem as the transition function $S^M(S_t,x_t,W_{t+1})$ can be broken into two parts: $S_t^x=S^{M,x}(S_t,x_t)$ and $S_{t+1}=S^{M,W}(S_t^x,W_{t+1})$. The resulting post-decision state-based VFA policy is given by
\begin{align}
\label{ApproxPolPDS}
X_t^{\pi}(S_t)=\argmax\limits_{x_t \in \mathcal{X}_t} \left(C(S_t,x_t)+\bar{V}^x_{t}(S^x_{t})\right),
\end{align}
where $\bar{V}^x_{t}(S^x_{t})$ serves as an approximation of $\mathbb{E}\left[V_{t+1}^*(S_{t+1})|S_t,x_t\right]$.

As mentioned previously, forward ADP is a far more common approach to fitting value function approximations for either $\bar{V}_{t+1}(S_{t+1})$ or $\bar{V}^x_{t}(S^x_{t})$, but the classical use of machine learning methods to approximate value functions with forward ADP has been found to work quite poorly in similar energy storage problems \citep{jiang2014comparison}. This paper instead utilizes a relatively seldom used ADP technique, backward ADP, as the resulting VFA-based policy achieves better performance.

Additional notation is necessary to describe the backward ADP algorithms in this paper. Let $\mathcal{S}_t$ and $\mathcal{S}_t^x$ be the set of all time $t$ pre- and post-decision states respectively, and $\hat{\mathcal{S}}_t^{\alpha}$ be a random sample of states in $\mathcal{S}_t$ sampled at rate $\alpha \in (0,1]$. Also, let $\mathbf{P}(S_{t+1}|S_t^x)$ be the transition probability for each pair $(S_t^x,S_{t+1})$ with $S^x_t \in \mathcal{S}_t^x$ and $S_{t+1} \in \mathcal{S}_{t+1}$. Finally, let $\mathcal{S}_{t+1}(S_t^x)$ be the set of all time $t+1$ pre-decision states such that $\mathbf{P}(S_t|S_t^x)>0$ (all the pre-decision states that can be reached from post-decision state $S_t^x$).

\subsection{Backward ADP with Lookup Tables}
Backward ADP algorithms resemble textbook backward dynamic programming, except that instead of looping over all states, only a sampled set of states are evaluated, and the results are used to create a value function approximation, replacing the exact lookup table value function used in classical discrete Markov decision processes. This reduces both the CPU time and memory necessary to compute and store value functions.

\begin{algorithm}[]
\scriptsize
\caption{Backward ADP with Post-Decision State Values Stored in Lookup Table Form}
\begin{algorithmic}
\STATE Choose sample rate $\alpha$.
\STATE Initialize terminal contributions $\bar{V}_T(S_T)$ for each $S_t \in \mathcal{S}_t$.
\STATE $\hat{\mathcal{S}}_T^{\alpha}=\left\lbrace\right\rbrace$
\FOR{each $S_{T-1}^x \in \mathcal{S}_{T-1}^x$}
\STATE Add $\left\lceil \alpha|\mathcal{S}_{T}(S_{T-1}^x)| \right\rceil$ states $S_T \in \mathcal{S}_{T}(S_{T-1}^x)$ to $\hat{\mathcal{S}}_T^{\alpha}$
\ENDFOR
\STATE Initialize $t \leftarrow T-1$.
\WHILE{$t\geq 0$}
\STATE{Expectation Step:}
\FOR{each $S_t^x \in \mathcal{S}_{t}^x$}
\STATE $\mathcal{S}_{t+1}^{sampled}=\mathcal{S}_{t+1}(S_{t}^x)\cap\hat{\mathcal{S}}_{t+1}^{\alpha}$
\STATE $p^{norm}=\sum\limits_{S_{t+1}\in \mathcal{S}_{t+1}^{sampled}} \mathbf{P}(S_{t+1}|S_t^x)$
\STATE $\bar{V}^x_t(S_t^x)=\frac{1}{p^{norm}}\sum\limits_{S_{t+1}\in \mathcal{S}_{t+1}^{sampled}} \mathbf{P}(S_{t+1}|S_t^x)\bar{V}_{t+1}(S_{t+1})$
\STATE Write $\bar{V}_{t}^x(S^x_t)$ to memory.
\ENDFOR
\STATE{Maximization Step:}
\IF{$t>0$}
\STATE $\hat{\mathcal{S}}_t^{\alpha}=\left\lbrace\right\rbrace$
\FOR{each $S_{t-1}^x \in \mathcal{S}_{t-1}^x$}
\STATE Add $\left\lceil \alpha|\mathcal{S}_{t}({S}_{t-1}^x)| \right\rceil$ states $S_t \in {\mathcal{S}}_{t}({S}_{t-1}^x)$ to $\hat{\mathcal{S}}_t^{\alpha}$
\ENDFOR
\FOR{$S_t \in \hat{\mathcal{S}}_t^{\alpha}$}
\STATE $\bar{V}_t(S_t)=\operatornamewithlimits{max}\limits_{x_t \in \mathcal{X}_t}\left\lbrace C(S_t,x_t)+\bar{V}_t^x(S_t^x)\right\rbrace$ where $S_t^x=S^{M,x}(S_t,x_t)$.
\ENDFOR
\ENDIF
\STATE $t \leftarrow t-1$
\ENDWHILE
\end{algorithmic}
\label{alg:Sampled MDP}
\end{algorithm}

First we sample the pre-decision state space to expedite the maximization step. Sampling pre-decision states at rate $\alpha$ forms the subset $\hat{\mathcal{S}}_{t+1}^\alpha$. Their values, $\bar{V}_{t+1}(S_{t+1})$ for $S_{t+1} \in \hat{\mathcal{S}}_{t+1}^\alpha$, are computed as usual by maximizing the one-step contribution plus the value of the resulting post-decision state over feasible decisions at each time step. 

The expectation step is also streamlined by computing post-decision state values $\bar{V}^x_{t}(S^x_{t})$ as a weighted average (with weights proportional to transition probabilities) of the values from the sampled set of pre-decision states. The approximations $\bar{V}^x_{t}(S^x_{t})$ are stored in a lookup table. This method, described in Algorithm \ref{alg:Sampled MDP}, requires careful sampling of the time $t+1$ pre-decision states such that each time $t$ post-decision state can reach at least one sampled next pre-decision state to form the approximation $\bar{V}^x_{t}(S^x_{t})$ for each $S^x_{t}$. This is accomplished by looping over each time $t$ post-decision state $S_{t}^x \in \mathcal{S}_t^x$ and sampling at minimum one time $t+1$ pre-decision state $S_{t+1}$ from $\mathcal{S}_{t+1}(S_{t}^x)$. 

The approach of storing post-decision state values in lookup table form is effective for problems in which the post-decision state space is much more compact than the pre-decision state space, as is the case in our energy storage problem. Of course, the accuracy of the approximate value functions depends on the sampling rate $\alpha$.

\subsection{Backward ADP using Parametric VFAs for Pre-Decision States}
\label{BADP_Par}

An alternative approach, laid out in Algorithm \ref{alg:Lin VFA}, is to again sample pre-decision states, but instead form a parametric approximation $\bar{V}_{t+1}(S_{t+1}|\theta_{t+1})$ given the set of sampled states and their resulting values following maximization over feasible actions. Assuming the chosen functional form $\bar{V}_{t+1}(S_{t+1}|\theta_{t+1})$ is parameterized by an $n$-dimensional vector $\theta_{t+1} \in \mathbb{R}^n$, if each $S_{t+1} \in \hat{\mathcal{S}}_{t+1}^{\alpha}$ has associated value $v(S_{t+1})$, we search for \begin{align}
\theta^*_{t+1}=\argmin\limits_{\theta_{t+1}} \sum\limits_{S_{t+1} \in \hat{\mathcal{S}}_{t+1}^{\alpha}} w(S_{t+1}) L\left(v(S_{t+1}),\bar{V}_{t+1}(S_{t+1}|\theta_{t+1})\right),
\end{align}
where $w(S_{t+1})$ are weights satisfying $\sum\limits_{S_{t+1} \in \hat{\mathcal{S}}_{t+1}^{\alpha}} w(S_{t+1})=1$ and $L(v(s),\bar{V}_{t+1}(s|\theta_{t+1}))$ is an appropriate loss function such as $L(v,\bar{v})=(v-\bar{v})^2$ ($\mathcal{L}_2$-norm) or $L(v,\bar{v})=|v-\bar{v}|$ ($\mathcal{L}_1$-norm). This chapter utilizes $L(v,\bar{v})=(v-\bar{v})^2$ and equally weights all samples. The approximation given by $\bar{V}_{t+1}(S_{t+1}|\theta^*_{t+1})$ is then used in the expectation step to approximate the values of all time $t+1$ pre-decision states. As a result, this approach does not require that at least one time $t+1$ pre-decision state $S_{t+1}$ from $\mathcal{S}_{t+1}(S_{t}^x)$ is sampled for each $S_{t}^x \in \mathcal{S}_t^x$ as in the lookup table representation. However, the practice of visiting each ${S}_{t-1}^x \in {\mathcal{S}}_{t-1}^x$ to produce at least one sample to add to $\mathcal{S}_{t}^{\alpha}$ is maintained here as it is an effective method of sampling different regions of the pre-decision state space given the compact nature of the post-decision state variable in the energy storage problems we explore (see Section \ref{ESP}). Additionally, the low-dimensional, best-fit parameter vector $\theta^*_{t+1}$ is stored in memory at each time step for later use when simulating the policy forward making this alternative much more memory efficient than a lookup table VFA. Note that this algorithm does still rely on a relatively compact post-decision state space (which exists in the storage problem we consider) as a lookup table representation, $\bar{V}^x_t({S}_t^x)$, is maintained for post-decision state values in the backward recursion.

\begin{algorithm}[]
\scriptsize
\caption{Backward ADP with Parametric VFAs}
\begin{algorithmic}
\STATE Choose sample rate $\alpha$ and VFA form $\bar{V}_{t}(S_t|\theta_t)$.
\STATE Initialize terminal contributions $\bar{V}_T(S_T)$ for each ${S}_T \in {\mathcal{S}}_T$.
\STATE $\hat{\mathcal{S}}_T^{\alpha}=\left\lbrace\right\rbrace$
\FOR{each $S_{T-1}^x \in \mathcal{S}_{T-1}^x$}
\STATE Add $\left\lceil \alpha|\mathcal{S}_{T}(S_{T-1}^x)| \right\rceil$ states $S_T \in \mathcal{S}_{T}(S_{T-1}^x)$ to $\hat{\mathcal{S}}_T^{\alpha}$
\ENDFOR
\STATE Initialize $t \leftarrow T-1$.
\WHILE{$t\geq 0$}
\STATE Solve $\theta^*_{t+1}=\argmin\limits_{\theta_{t+1}} \sum\limits_{S_{t+1} \in \hat{\mathcal{S}}_{t+1}^{\alpha}} w(S_{t+1}) L\left(v(S_{t+1}),\bar{V}_{t+1}(S_{t+1}|\theta_{t+1})\right)$
\STATE Write $\theta^*_{t+1}$ to memory.
\STATE Expectation Step:
\FOR{each ${S}_t^x \in {\mathcal{S}}_{t}^x$}
\STATE $\bar{V}^x_t({S}_t^x)=\sum\limits_{S_{t+1}\in {\mathcal{S}}_{t+1}} \mathbf{P}(S_{t+1}|{S}_t^x)\bar{V}_{t+1}(S_{t+1}|\theta_{t+1}^*)$.
\ENDFOR
\STATE Maximization Step:
\IF{$t>0$}
\STATE $\hat{\mathcal{S}}_t^{\alpha}=\left\lbrace\right\rbrace$
\FOR{each $S_{t-1}^x \in \mathcal{S}_{t-1}^x$}
\STATE Add $\left\lceil \alpha|\mathcal{S}_{t}({S}_{t-1}^x)| \right\rceil$ states $S_t \in {\mathcal{S}}_{t}({S}_{t-1}^x)$ to $\hat{\mathcal{S}}_t^{\alpha}$
\ENDFOR
\FOR{${S}_t \in \hat{\mathcal{S}}_t^{\alpha}$}
\STATE $v({S}_t)=\operatornamewithlimits{max}\limits_{x_t \in \mathcal{X}_t}\left\lbrace C({S}_t,x_t)+\bar{V}^x_t({S}_t^x)\right\rbrace$ where ${S}_t^x=S^{M,x}({S}_t,x_t)$.
\ENDFOR
\ENDIF
\STATE $t \leftarrow t-1$
\ENDWHILE
\end{algorithmic}
\label{alg:Lin VFA}
\end{algorithm}

A commonly used architecture for VFAs is a linear model (that is, linear in the parameters), $\bar{V}_t(S_t|\theta_t)=\sum\limits_{k=0}^n \theta_{t,k}\phi_k(S_t)$ where $\phi_0(S_t)=1$ and $\phi_k(S_t)$ for $k \in \left\lbrace 1, 2, ..., n \right\rbrace$ are basis functions of the states that must be chosen beforehand and should be appropriate for the application. Picking these is not a trivial task and if basis functions are chosen poorly, solution quality will suffer. This paper utilizes basis functions that are relevant to the energy storage problem which are described in greater detail in Section \ref{Numerical Results}.

Note that, despite only utilizing a linear VFA in this paper, the structure of Algorithm \ref{alg:Lin VFA} can be generalized to work for many parametric or even nonparametric approximations for $\bar{V}_t(S_t)$. This requires using the appropriate method for fitting the function given the pairs of sampled pre-decision states and their values and storing any parameters or data necessary for simulating forward.

\subsection{Backward ADP for the Energy Storage Problem}
\label{ESP}

By using the hidden semi-Markov crossing state model presented in \cite{durante2017}, we create a more realistic energy storage problem, but in the process form a partially observable Markov decision problem (POMDP) as the crossing states are partially hidden. It is possible to find an optimal policy for a POMDP by incorporating the belief state in the state variable \citep[see][]{KAELBLING199899, murphy2000survey}. However, it is often computationally difficult (if not impossible) to solve belief MDPs exactly as the dimensionality of the problem is often greatly increased. Policies based on finding value functions for the underlying completely observable MDP are simple, yet effective, heuristic approaches to circumventing this issue \citep[see, for example,][]{cassandra1998exact}. This is the approach we use here and the ``optimal policies" in section \ref{Numerical Results} refer to VFAs that are optimal with respect to an underlying fully observable MDP.

Recall that when simulating forward, we have belief states $B_t^{E}$ and $B_t^{P}$, giving the operator's distribution of belief about the crossing states given the history of the process up to time $t$. This information is not available in the backward pass. Instead, value functions are fit to each possible resource state conditioned on each information state in the backward pass. The value of being at resource level $R_t$ at time $t$ given the wind and price processes are in information states $I_t^E$ and $I_t^P$ is given by
\begin{align}
\bar{V}_t(R_t,I_t^E,I_t^P)=\max\limits_{x_t \in \mathcal{X}_t} \left(C(R_t,I_t^E,I_t^P,x_t)+\bar{V}^x_{t}(R^x_t,I_t^E,I_t^P)\right),
\end{align} where $\bar{V}^x_{t}(R^x_t,I_t^E,I_t^P)=\mathbb{E}\left[\bar{V}_{t+1}(R_{t+1},I_{t+1}^E,I_{t+1}^P)|R_t^x,I^E_t,I^P_t \right]$. However, as mentioned in Section \ref{CompactInfo}, finding values for all possible system states using the full information state is very computationally expensive.

These computational issues are addressed by utilizing the compact information states introduced in Section \ref{CompactInfo} in combination with backward ADP. We need to introduce additional notation to handle this. Let $\tilde{S}_t=(R_t,E_t,P_t,\tilde{I}_t^E,\tilde{I}_t^P)$ and $\tilde{S}_t^x=(R_t^x,\tilde{I}_t^E,\tilde{I}_t^P)$ be compact time $t$ pre- and post-decision states respectively using the compact information states for the wind and price processes $\tilde{I}^{E}_t$ and $\tilde{I}^{P}_t$ from Section \ref{CompactInfo}. The sets of all possible compact time $t$ pre- and post-decision states are then represented by $\tilde{\mathcal{S}}_t$ and $\tilde{\mathcal{S}}^x_t$. Using these compact state variables, the post-decision to next pre-decision state transition matrix is defined as \begin{align}
\mathbf{P}(\tilde{S}_{t+1}|\tilde{S}_t^x)=\mathbf{1}_{\left\lbrace R_{t+1}=R_t^x \right\rbrace} \mathbf{P}(\tilde{I}_{t+1}^E,E_{t+1}|\tilde{I}_t^E)\mathbf{P}(\tilde{I}_{t+1}^P,P_{t+1}|\tilde{I}_{t}^P),
\end{align} where, noting $\hat{E}_{t+1}=E_{t+1}-f_{t+1}^E$,

\begin{align}
\mathbf{P}(\tilde{I}_{t+1}^E=i',E_{t+1}|\tilde{I}_t^E=i)=\mathbf{1}_{\left\lbrace \hat{E}_{t+1} \in \hat{E}_{t+1}^b \right\rbrace}\tilde{\mathbb{P}}(i'|i) \times
&\begin{cases}
\mathbf{P}(\hat{E}_{t+1}|i,\hat{E}^b_t) &\text{if } i'=i \\
\mathbf{P}(\hat{E}_{t+1}|i',t+1 \in \mathcal{C}^{E,U} \cup \mathcal{C}^{E,D})&\text{otherwise}, 
\end{cases}
\end{align} with matrix $\tilde{\mathbb{P}}(i'|i)$ defined in Section \ref{CompactInfo}. $\mathbf{P}(\tilde{I}_{t+1}^P,P_{t+1}|\tilde{I}_t^P)$ is defined similarly.

Replacing $S_t,S_t^x,\mathcal{S}_t$ and $\mathcal{S}^x_t$ with $\tilde{S}_t, \tilde{S}^x_t, \tilde{\mathcal{S}}_t$, and $\tilde{\mathcal{S}}^x_t$ in Algorithms \ref{alg:Sampled MDP} and \ref{alg:Lin VFA}, we find VFAs for the simplified MDP. Despite using these simplified state variables in the backward pass, we utilize the full state variables to make decisions moving forward in time according to the policy given by Equation \eqref{ApproxPol}. If Algorithm \ref{alg:Sampled MDP} was used, the expectation $\mathbb{E}[\bar{V}_{t+1}(S_{t+1})|S_t,x_t]$ in this policy is computed based on approximate values for simplified post-decision states: $\bar{V}^x_t(\tilde{S}_t^x)$. Note that $\bar{V}^x_t(\tilde{S}_t^x)=\bar{V}^x_t(R_t^x,\tilde{I}_t^E,\tilde{I}_t^P)=\bar{V}^x_t(R_t^x,(I_t^{C,E},\hat{E}_t^b),(I_t^{C,P},\hat{P}_t^b))$ is an approximation of the expected value of the downstream pre-decision state conditioned on the crossing state variables $I_t^{C,E}$ and $I_t^{C,P}$. By applying the law of total expectation using our time $t$ belief states, we have
\begin{align}
\label{ESPDS}
\mathbb{E}[\bar{V}_{t+1}(S_{t+1})|S_t,x_t]=&\sum\limits_{i\in \mathcal{I}^{C,E}} \mathbf{P}(I_t^{C,E}=i)\sum\limits_{j\in \mathcal{I}^{C,P}}\mathbf{P}(I_t^{C,P}=j)\left(\bar{V}_t^x(R_{t}^x,(i, \hat{E}_t^b),(j,\hat{P}_t^b))\mathbf{1}_{\left\lbrace R_{t}^x = S^{M,x}(R_t,x_t) \right\rbrace}\right),
\end{align} where $\mathbf{P}(I_t^{C,E}=i)$ and $\mathbf{P}(I_t^{C,P}=j)$ (contained in $B_t^E$ and $B_t^P$) give the probability that the stochastic wind and price processes are in each crossing state at time $t$ based on their history. If, instead, parametric VFAs were calculated and stored using Algorithm \ref{alg:Lin VFA}, the policy can be computed by simply replacing $\bar{V}^x_t(R_t^x,(I_t^{C,E},\hat{E}_t^b),(I_t^{C,P},\hat{P}_t^b))$ in Equation \eqref{ESPDS} with $\mathbb{E}[\bar{V}_{t+1}(S_{t+1}|\theta_{t+1}^*)|R_t^x,(I_t^{C,E},\hat{E}_t^b),(I_t^{C,P},\hat{P}_t^b)]$.

\section{Numerical Results}
\label{Numerical Results}

In this section we observe the trade-off between computation time and performance for the backward ADP methods at different sampling rates. We show that, in comparison to an exact solution to the simplified MDP, we can reduce computation time significantly via backward ADP without much degradation in solution quality. Additionally, we show that backward ADP outperforms a common forward ADP method, approximate policy iteration, a parametric policy function approximation, and direct lookahead policies in this problem setting.

The above claims are substantiated by the results from several test cases, for which input parameters are varied. We test different combinations of battery sizes, battery charge rates (where we assume $\rho=\rho^{ch}=\rho^{dch}$), wind power forecasts, temperature forecasts, and load profiles. Table \ref{table:Test Case Table} describes each test case. Also displayed is the C-rate of the battery, where C$=(\rho \times \text{number of time periods per hour})/R^{max}$, which is the maximum discharge rate relative to battery capacity. Common C-rates range from $0.1$C (slow charging) to $2$C and higher (fast).

\begin{table}[h]
\TABLE
{Test cases for the energy storage problem. The temperature forecast (Hot, Avg, Cool) also determines the load profile. The numbers in the $\left\lbrace f_t^E\right\rbrace_{t=0}^T$ row refer to the wind forecasts shown in Figure \ref{CrossingTimePic}. \label{table:Test Case Table}}
{\begin{tabular}{c || c c c c c c c c c c c c} % centered columns (4 columns)
\hline
\hline
 %inserts double horizontal lines
Case & 1 &2 &3&4&5&6&7&8&9&10&11&12  \\ [0.5ex] % inserts table
%heading
\hline % inserts single horizontal line
Temperature & Avg &Cool& Hot& Avg &  Cool & Cool &Hot & Hot & Avg & Avg &Hot&Hot\\
$\left\lbrace f_t^E\right\rbrace_{t=0}^T$ & 1&1&1&1&2&2&1&1&2&2&2&2 \\
$R^{max}$ (MWh)  & 5.00& 5.00 & 5.00 &5.00& 4.00&4.00&3.33&3.33&2.00&4.00&5.00&2.50\\
$\rho$ (kWh)& 83.3&83.3&83.3&416.7&333.3&83.3&416.7&83.3&41.7&83.3&83.3&83.3\\
C-rate & 0.20C&0.20C&0.20C&1.00C&1.00C&0.25C&1.50C&0.30C&0.25C&0.25C&0.20C&0.40C \\
\hline
\end{tabular}}
 % is used to refer this table in the text
{}
\end{table}

Across the test cases, we set $m^E=3$ in the wind model which determines the number of crossing time duration bins ($D_t^{E}$). Additionally, we let $n^E=3$, fixing the number of error bins ($\hat{E}_t^{b}$) for each crossing state. Wind power ($E_t$), which can take on continuous values, is discretized to form a discrete MDP. Using a uniform interval of $100$ kW, we allow $51$ possible values for $E_t$, evenly spaced from $0$ MW to $E^{max}=5$ MW. Similarly, we let $m^P=1$ and $n^P=4$ in the price model and discretize $P_t$ evenly using an interval of $\$2$/MWh from $P^{min}$ to $P^{max}$ where these represent the lowest and highest prices observed in training data. Here, $P^{min}=-\$312$/MWh and $P^{max}=\$780$/MWh, resulting in $547$ possible values $P_t$ may take on. The state of charge in the battery, $R_t$, is also discretized evenly from $0$ MWh to $R^{max}$. The number of charge states will, in general, depend on the charge rate, $R^{max}$, and the program time step, though the time step is fixed at five minutes in these problems. The number of battery charge states is between 30 and 60 in the test cases (slower charging, larger batteries need more states). Additionally, $\eta$ is set to 1.0, and we assume the battery starts out with a 50 percent charge. Finally, each test case has an optimization horizon of $T=288$ as we make an energy allocation decision every five minutes over the course of one day.

\subsection{Policies Tested}
Backward ADP methods are used to find value functions in each test case with various sampling rates $\alpha$. Type Lookup-$\alpha$ is the approach described in Algorithm \ref{alg:Sampled MDP} which uses a lookup table representation of the value function. Type Lin-$\alpha$ represents the method in Algorithm \ref{alg:Lin VFA} using linear parametric VFAs and basis functions given by: $\phi_1(S_t)=(C_t^E-0.5)(D_t^E+1)$, $\phi_2(S_t)=E_t$, $\phi_3(S_t)=E^2_t$, $\phi_4(S_t)=(C_t^P-0.5)(D_t^P+1)$, $\phi_5(S_t)=P_t$, $\phi_6(S_t)=P^2_t$, $\phi_7(S_t)=R_t$, $\phi_8(S_t)=R_t^2$, $\phi_9(S_t)=E_t P_t$, $\phi_{10}(S_t)=E_t R_t$, and $\phi_{11}(S_t)=P_t R_t$, where $\phi_1(S_t)$ and $\phi_4(S_t)$ are chosen such that the longer down-crossing states evaluate to lower numbers and the longer up-crossing states evaluate to higher numbers.

%This mimics the general relationship between the wind process crossing states and their values: the longer the down-crossing duration bin, the smaller the value (as there is less wind energy for longer), while up-crossing states increase in value as the duration bin increases.
%
%This relationship is not as clear for the price process. A high price may be beneficial at some times, such as periods when there is an excess of wind and the battery is well charged. In the opposite scenario, when there is little wind and an empty battery, higher prices will result in larger costs. Though, this is not a problem as there is a best fit $\theta_t^*$ for each time period. $\theta_{t,4}^*$ will determine whether the correlation between price crossing state and state value is positive or negative at each point in time.

The backward ADP policies are also compared to several other policies. One is the optimal policy using the value functions resulting from exact backward dynamic programming. Another is approximate policy iteration (API) utilizing linear regression to fit value functions to post-decision states. This is a form of forward ADP (for a more thorough discussion of API, see \cite{jiang2014comparison}). The performance of this policy will improve over time as the number of training iterations increases. We restrict training time to 12 hours, often much longer than the backward ADP algorithms take to run in the test problems, and use the resulting policy at the end of this training period. It is also worth noting that, as with backward ADP with linear VFAs, policy performance is highly dependent on the selection of proper basis functions here as well.

Another policy tested is a carefully tuned buy-low, sell-high policy function approximation (PFA) defined by $\theta=(\theta^H,\theta^L)$ with $\theta^H>\theta^L$, given by \begin{equation}
\begin{aligned}
X_t^{PFA}(S_t|\theta)=\begin{cases}
x_t^{EL}=\min\left\lbrace{L_t,E_t}\right\rbrace\\
x_t^{RL}=\begin{cases}\min\left\lbrace{L_t-x_t^{EL},\min\left\lbrace R_t,\rho^{dch}\right\rbrace}\right\rbrace &\text{ if } P_t>\theta^H\\
0 &\text{ if } P_t \leq \theta^H\\
\end{cases}\\
x_t^{GL}=L_t-x_t^{EL}-x_t^{RL}\\
x_t^{ER}=\min\left\lbrace{E_t-x_t^{EL},\rho^{ch},R^{max}-R_t}\right\rbrace\\
x_t^{GR}=\begin{cases} \min\left\lbrace\rho^{ch}-x_t^{ER},R^{max}-R_t-x_t^{ER}\right\rbrace &\text{ if } P_t<\theta^L\\
0 &\text{ if } P_t \geq \theta^L\\
\end{cases}\\
x_t^{RG}=\begin{cases}\min\left\lbrace{R_t-x_t^{RL}, \rho^{dch}-x^{RL}}\right\rbrace &\text{ if } P_t>\theta^H\\
0 &\text{ if } P_t \leq \theta^H,\\
\end{cases}\\
\end{cases}
\end{aligned}
\end{equation} where $\theta$ is manually tuned via grid search. The $\theta$ resulting in the highest average objective function value in each case is chosen.

Finally, we compare to direct lookahead policies, also referred to as rolling or receding horizon procedures, of varying horizons and forecast accuracy. For a lookahead horizon $H$, the DLA policy is given by
\begin{equation}
X_t^{DLA}(S_t|H)=\argmax\limits_{x_t \in \mathcal{X}_t}\left(C(S_t,x_t)+\max\limits_{\tilde{x}_{t,t+1},...,\tilde{x}_{t,t+H}}\sum\limits_{t'=t+1}^{t+H}C(S_{tt'},\tilde{x}_{tt'})\right),
\end{equation}where $R_{t,t+1}=R_{t}+\eta(x_t^{GR}+x_{t}^{ER})-x_{t}^{RG}-x_{t}^{RL}$, $C(S_{tt'},\tilde{x}_{tt'})=f^P_{tt'}(L_{t'}-\tilde{x}_{tt'}^{GR}-\tilde{x}_{tt'}^{GL}+\eta \tilde{x}_{tt'}^{RG})$, and lookahead decisions $\tilde{x}_{tt'}$ for $t'=t+1,...,t+H$ are subject to
\begin{samepage}
\begin{align}
\tilde{x}_{tt'}^{EL}+\tilde{x}_{tt'}^{ER} &\leq f^{E}_{tt'},\\
\tilde{x}_{tt'}^{GL}+\tilde{x}_{tt'}^{EL}+\eta \tilde{x}_{tt'}^{RL} &= L_{t'},\\
\tilde{x}_{tt'}^{RG}+\tilde{x}_{tt'}^{RL} &\leq \min(R_{tt'},\rho^{dch}),\\
\tilde{x}_{tt'}^{ER}+\tilde{x}_{tt'}^{GR} &\leq \min(\rho^{ch},R^{max}-R_{tt'}),\\
R_{t,t'+1}&=R_{tt'}+\eta(\tilde{x}_{tt'}^{GR}+\tilde{x}_{tt'}^{ER})-\tilde{x}_{tt'}^{RG}-\tilde{x}_{tt'}^{RL},\\
\tilde{x}_{tt'}^{GL},\tilde{x}_{tt'}^{GR},\tilde{x}_{tt'}^{RG},\tilde{x}_{tt'}^{EL},\tilde{x}_{tt'}^{ER},\tilde{x}_{tt'}^{RL} &\geq 0.
\end{align}
\end{samepage}Despite optimizing over a longer horizon, the policy implements only the decision $x_{t}$ at time $t$ and, moving forward in time, re-optimizes at each time step by solving a new linear program.

In practice, updated forecasts at time $t$ for price, $f_{tt'}^P$, and wind, $f_{tt'}^E$, over the horizon $t'=t+1,...,t+H$ would be used in the lookahead model. Here, we instead contrive a forecast based on perfect knowledge of the future up to the end of the horizon $H$. A length $H$ zero-mean Gaussian random vector $\epsilon^E$ with covariance matrix $\Sigma^E$ is added to the actual realization of wind, with
\begin{align}
\Sigma^E=\begin{bmatrix}
(\sigma_1^E)^2 & \gamma^E \sigma^E_1 \sigma^E_2 & (\gamma^E)^2 \sigma^E_1 \sigma^E_3 & \dots & (\gamma^E)^{H-1} \sigma^E_1 \sigma^E_H\\
\gamma^E \sigma^E_2 \sigma^E_1 & (\sigma_2^E)^2 & \gamma^E \sigma^E_2 \sigma^E_3 & \dots & (\gamma^E)^{H-2} \sigma^E_2 \sigma^E_H\\
(\gamma^E)^2 \sigma^E_3 \sigma^E_1 & \gamma^E \sigma^E_3 \sigma^E_2 & (\sigma^E_3)^2 & \dots & (\gamma^E)^{H-3} \sigma^E_3 \sigma^E_H\\
\vdots& \vdots &\vdots &\ddots &\vdots\\
(\gamma^E)^{H-1} \sigma^E_H \sigma^E_1 & (\gamma^E)^{H-2} \sigma^E_H \sigma^E_2 & (\gamma^E)^{H-3} \sigma^E_H \sigma^E_3 & \dots &(\sigma^E_H)^2
\end{bmatrix},
\end{align}where $\sigma^E_{i}$ is the standard deviation of the simulated error for lead time $i$ and $\gamma^E=\in [0,1]$ is the correlation between successive error terms (the lag-1 autocorrelation) in the training set of forecast errors. Letting $f^E_{tt'}=E_{t'}+\epsilon^E_{t'-t}$, we see that the actual realization of wind can be viewed as a deviation from the contrived forecast, with $\mathbb{E}[f^E_{tt'}]=E_{t'}$. As error variance increases with forecast lead time, we ensure $\sigma^E_1<\sigma^E_2<\dots<\sigma^E_H$. This is accomplished by letting $\sigma_{i}^E=i c^E s^E$ where $s^E$ is the observed standard deviation of forecast errors from training data and $c^E$ is a constant which allows us to control the accuracy of the constructed forecast. The same procedure is carried out to construct artificial price forecasts where the superscript $P$ takes the place of the superscript $E$ for all the relevant variables and constants. A DLA policy with lookahead horizon $H$ and $c=c^E=c^P$ is referred to as DLA-$H$-$c$ where a smaller $c$ results in more accurate forecasts.

\subsection{Policy Performance Results}

The policies described in the previous section are simulated using a set of 100 sample paths of wind energy forecast errors and LMPs generated from the stochastic base model (the HSMMs). Table \ref{table:Performance Table} reports each policy's performance as the cumulative contribution over 100 trials as a percent of the cumulative contribution earned using the exact solution to the discretized MDP.

\begin{table}[h]
\TABLE
{Mean performance of policies in the various test cases over 100 trials, reported as a percent of the optimal policy. Shown on the far right is the average time in hours needed to compute value functions for each algorithm (or time allotted for policy search), but note that this is highly problem-dependent and can vary greatly between test cases. The optimal policy took an average of \textit{11.3} hours to compute.
\label{table:Performance Table}}
 % used for centering table
{\begin{tabular}{c || c c c c c c c c c c c c | c || c} % centered columns (4 columns)
\hline\hline %inserts double horizontal lines
Case&1&2&3&4&5&6&7&8&9&10&11&12&Avg&CPU (hrs)\\ [0.5ex] % inserts table
%heading
\hline % inserts single horizontal line
Lookup-.01&99.1&96.7&98.0&88.8&98.0&100.0&93.2&98.4&99.6&99.4&98.8&98.9&97.4&0.41\\
Lin-.01&96.2&96.4&96.9&91.8&88.2&95.2&95.4&98.9&98.9&98.2&98.3&99.5&96.2&1.62\\
Lookup-.10&100.1&99.5&99.3&97.1&99.7&100.2&97.2&99.4&100.0&100.1&99.6&99.8&99.3&0.67\\
Lin-.10&96.3&96.5&98.1&91.1&88.3&95.2&94.9&98.9&98.9&98.2&99.0&99.2&96.2&2.72\\
API&82.7&77.9&79.5&57.3&76.2&90.8&50.5&80.4&94.7&90.0&86.4&86.9&79.5&12.0\\
PFA&93.6&92.3&93.4&71.4&80.6&91.3&72.3&93.4&97.3&95.2&96.0&94.6&89.3&5.14\\
DLA-24-.05&87.7& 86.6 & 87.4&73.5&84.8&93.3&90.4&89.9&96.8&92.2&92.3&95.3&89.2&N/A\\
DLA-72-.01&91.1&91.5& 90.8&86.1&99.4&95.5&101.7&93.0&97.9&94.4&94.2&94.3&94.2&N/A\\
\end{tabular}}
{}
\end{table}

From Table \ref{table:Performance Table}, we see that both backward ADP methods frequently outperform approximate policy iteration, a policy function approximation, and the deterministic lookahead policies (even when the forecasts are quite accurate). Also observe that the lookup table method performs better than the linear VFA at both sampling rates in most cases. This may be because the linear structure cannot capture the complex relationships between states and their values that exist in this problem, even with a large set of sampled states. The linear VFA method is likely to perform better in a problem where the value functions exhibit more structure, such as a known polynomial form. However, note that the performance of the linear VFA is approximately the same for both the larger and smaller sampling rate $\alpha$, whereas this is not the case for the lookup table method as performance degrades with decreasing $\alpha$. This is an encouraging sign for the use of Algorithm \ref{alg:Lin VFA} and linear VFAs in MDPs with much larger pre-decision state spaces where Algorithm \ref{alg:Sampled MDP} may not work well as it may be necessary to sample at an even smaller rate $\alpha$.

Table \ref{table:Performance Table} also reports the average CPU time in hours required to compute value functions for each dynamic programming algorithm (in the API case it is the maximum allotted training time). In general, sampling in the backward ADP algorithms will result in much quicker CPU times (as expected) without sacrificing much performance in comparison to the full MDP solution. This is useful as the backward pass can be performed closer to the time at which it is used to control the system. The accuracy of the input load, wind, and temperature forecasts will be improved as a result of the shorter forecast lead times.

Finally, we note that the storage problems considered in this paper required between $.075$ GB and $.15$ GB of memory to store value functions in lookup table form for post-decision states. This is a significant improvement over storing pre-decision state values in lookup table form; however this number can also grow quickly as the dimensionality of the post-decision state increases. Thus, if memory storage becomes an issue, the linear VFA form may be preferable as it only requires that a low- dimensional parameter vector be stored at each time step, requiring something on the order of $10$ KB to $100$ KB of storage space depending on the dimension of the parameter vector (regardless of the dimensionality of the state space).

%The advantage over the deterministic lookaheads is especially pronounced if the variance of the errors from the rolling forecasts is high and increasing quickly with lead time. In a rolling horizon DLA, batteries with faster charge rates can quickly compensate for poor decisions made in the past by making good decisions in the present based on relatively accurate forecasts of the very near future (for example charging right before a price spike is anticipated in the near future). The slow charging battery, on the other hand, will not be able to do so and must plan far in advance for the possibility of price spikes and low renewable power production.

%It is clear that the backward ADP policies are most useful for systems with slower charging batteries, performing quite close to the posterior optimal solution in these cases. 

%Note that we have not considered the possibility of re-optimizing at finer time scales using backward ADP methods with periodic updates of the forecasts. Given we want to compute new value functions every $k$ time periods, we can choose a sampling rate, backward ADP method, and optimization horizon that would make this computation possible in the allotted time. This hybrid DLA-VFA policy has the potential to produce high quality solutions if done correctly.

\section{The Value of Hidden Semi-Markov Models}
\label{Model Selection}

Finally, in this section we observe that the choice of stochastic model has an effect on the quality of the solution. Specifically, we present results supporting the claim that explicitly modeling crossing time behavior leads to the development of more robust policies for storage systems. We show this by training VFAs on two different models for wind power forecast errors that take into account intertemporal correlation, one that replicates crossing times well and one that does not, and then evaluating the resulting policies on sample paths drawn from history. The models used to train VFAs are:
\begin{henumerate}
\item{A first-order autoregressive (AR) model: $\hat{E}_t=\gamma^E \hat{E}_{t-1}+\mathcal{N}(0,(s^{E,resid})^2)$ where $\gamma^E$ is the lag-1 coefficient and $s^{E,resid}$ is the standard deviation of the residuals when fitting this model to training data. This model does not replicate crossing times well (see Figure \ref{WindCrossingTime}).}
\item{The crossing state HSMM described in this paper, which replicates crossing times accurately.}
\end{henumerate} The AR model is restricted to a first order model (an AR(1) model) to keep the state variable compact enough to perform backward dynamic programming. This is a fairly common assumption seen in the literature when dynamic programming approaches are considered. For this comparison, the price process is modeled with the crossing state model described in Section \ref{WPM} (with more detail in the online Appendix) in both cases.

We train VFAs using the two wind models in three cases from Table \ref{table:Test Case Table}. Additionally, for each test case we train the wind power model on forecast error data sets from two different Great Plains wind farms from different months of the year (see Section \ref{WPM} for a description of the data). Value functions are then fit to post-decision states in lookup table form using either exact or approximate backward dynamic programming for each model. For the HSMM, we use Algorithm \ref{alg:Sampled MDP} with a sampling rate of $\alpha=0.10$. The lookup table VFAs for the AR(1) error model are found using exact backward dynamic programming. Thus, the resulting policy is an optimal policy under this modeling assumption. Note that the full solution for the HSMM could be computed in about the same amount of time as the full solution using the AR model. However, the purpose of this test is to show that the benefits of using a better stochastic model can, depending on the application, outweigh the negative effects of approximations in the optimization algorithm that may be necessary to efficiently compute a solution. Thus, the approximate solution is used.

In each case, the policies are then evaluated on two sets of 25 wind power and LMP sample paths: one \textit{Typical} set and one \textit{Worst Case} set. The \textit{Typical} sample paths are a set of 25 historical wind power forecast error and electricity price sample paths chosen at random. Let this set be $\Omega^{Typ}$ and each individual sample path be $\omega^{Typ}\in \Omega^{Typ}$. Testing is then repeated on a set of 25 \textit{Worst Case} wind and electricity price sample paths created by pairing a wind forecast error sample path from one of the five days in history with the lowest average forecast error (most negative) with an LMP path from one of the five days with the highest average LMP observed in history. This produces conditions that could lead to especially high operating costs. Let this set be $\Omega^{WC}$ and each individual sample path be $\omega^{WC}\in \Omega^{WC}$. Let $F^{\pi,Typ}$ and $F^{\pi,WC}$ be the mean additional profit earned by following policy $\pi$ outside of the baseline profit earned from satisfying the demand (as any feasible policy would receive this portion of the profit) in both \textit{Typical} and \textit{Worst Case} scenarios. Letting $\tilde{C}(S_t,x_t)=P_t(\eta x_t^{RG}-x_t^{GR}-x_t^{GL})$ be the contribution function with the baseline profit subtracted, $F^{\pi,Typ}$ is calculated as
\begin{equation}
\label{PolVal}
F^{\pi,Typ}=\frac{1}{|\Omega^{Typ}|}\sum\limits_{\omega^{Typ}\in \Omega^{Typ}}\sum\limits_{t=0}^{T} \tilde{C}\left(S_t(\omega^{Typ}),X_t^{\pi}(\omega^{Typ})\right),
\end{equation} where $S_{t+1}(\omega^{Typ})=S^M\left(S_t(\omega^{Typ}),X^{\pi}_t(S_t(\omega^{Typ})),W_{t+1}(\omega^{Typ})\right)$ (see Section \ref{Transition Function} for the definition of the transition function). $F^{\pi,WC}$ is also calculated using Equation \eqref{PolVal}, but with all ``Typ" superscripts replaced with ``WC." Table \ref{table:Model Table} reports $F^{\pi,Typ}$ and $F^{\pi,WC}$ in each test case for both the AR(1)-trained and HSMM-trained policies.

\begin{table}
\TABLE
{In several test cases, and for two different wind farms, value functions are trained assuming either an AR(1) model for wind power forecast errors or the crossing state HSMM. For each test case, $F^{\pi,Typ}$ and $F^{\pi,WC}$, the mean baseline-shifted profits in both \textit{Typical} and \textit{Worst Case} scenarios (calculated using Equation \eqref{PolVal}), are reported in dollars when the resulting policies are evaluated on each set of sample paths. Bold indicates better performance in each case.
\label{table:Model Table}}
 % used for centering table
{\begin{tabular}{c || c c | c c | c c | c c | c c | c c } % centered columns (4 columns)
\hline\hline 
\multirow{3}{*}{Wind Power Model} & \multicolumn{4}{c}{Test Case 1} & \multicolumn{4}{c}{Test Case 5} & \multicolumn{4}{c}{Test Case 9} \\
\cline{2-13}
& \multicolumn{2}{c}{Farm 1} & \multicolumn{2}{c}{Farm 2} & \multicolumn{2}{c}{Farm 1} & \multicolumn{2}{c}{Farm 2} & \multicolumn{2}{c}{Farm 1} & \multicolumn{2}{c}{Farm 2} \\
\cline{2-13}
 & Typ & WC & Typ & WC & Typ & WC & Typ & WC & Typ & WC & Typ & WC \\
\hline % inserts single horizontal line
 AR(1) & -960 & -2820 & -399 & -2143 & 55 & -306 & \bf{581} & -568 &-425 &-1325&\bf{-228}&-1989 \\
  HSMM & \bf{-959} & \bf{-2769} & \bf{-378} & \bf{-2059} & \bf{79} & \bf{-41} & 342 & \bf{-392} & \bf{-420} & \bf{-1300}& -242&\bf{-1902}\\

 % is used to refer this table in the text
\end{tabular}}
{}
\end{table}

\begin{figure}
\centering
\includegraphics[width=\columnwidth, height=2.8 in]{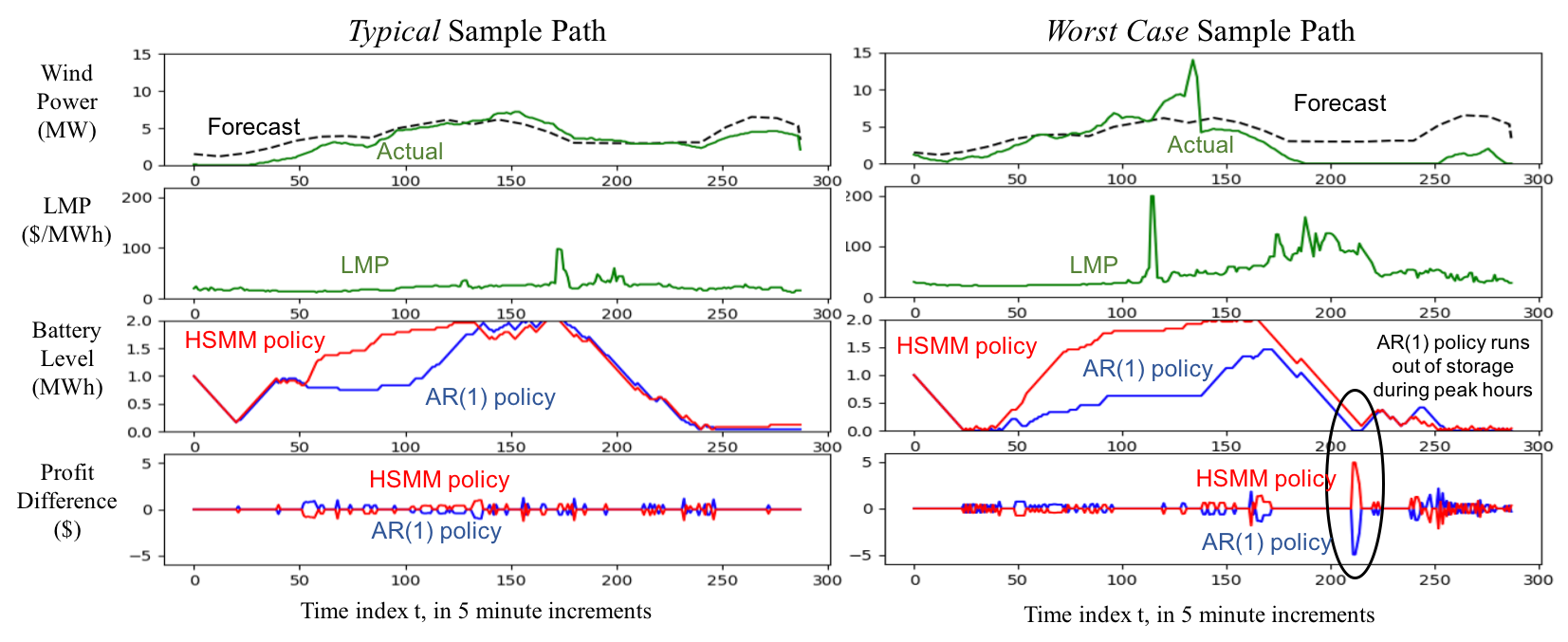}
\caption{For one \textit{Typical} (left) and one \textit{Worst Case} (right) sample path in test case 9, the following plots are shown, from top to bottom: 1) the forecasted and actual wind power output, where wind power falls to 0 MW during peak hours in the \textit{Worst Case} path, 2) the LMP path, which is higher and exhibits more spikes in the \textit{Worst Case} path, 3) the battery level at each point in time following both policies, where the HSMM-trained policy plans extra storage early in the day and does not run out in the \textit{Worst Case} scenario as the AR(1)-trained policy does (see the circled region), and 4) the deviation from the mean profit at time $t$ earned by both policies, where profit differences tend to even out in the \textit{Typical} path, but the robust HSMM-based policy avoids high costs during peak hours in the \textit{Worst Case} sample path.}
\label{Profits}
\end{figure} 

Observe from Table \ref{table:Model Table} that while the AR(1)-trained policy may perform similarly (or perhaps better) in \textit{Typical} scenarios, the crossing state HSMM-trained policy performs better across the test cases in \textit{Worst Case} scenarios. Also note that the profit difference is quite significant in \textit{Worst Case} scenarios. This indicates that training on the HSMM results in more robust policies. Insight as to why the HSMM produces a robust policy can be garnered from Figure \ref{Profits}, in which the AR(1)-trained and HSMM-trained policies are simulated for a single \textit{Typical} and \textit{Worst Case} sample path. For both cases, each policy's deviation from the average profit at time $t$, along with the corresponding battery level are plotted. The AR(1) model tends to underestimate the amount of time wind stays below its forecast (see Figure \ref{WindCrossingTime}) and thus does not emphasize storage for \textit{Worst Case} scenarios. Thus, we see the AR(1)-trained policy tending to sell more energy early in the day, expecting enough wind to be available in later periods. The HSMM-based policy, which replenishes storage levels to a greater extent in preparation for the possibility of extended periods with little to no renewable power, will perform worse early in the day, but profit much more during the peak hours when a \textit{Worst Case} scenario occurs. This is evident in the circled region of Figure \ref{Profits} which highlights a period when the AR(1)-trained policy runs out of battery storage during peak hours in the \textit{Worst Case} scenario and suffers high operating costs while the HSMM-based policy does not. Meanwhile, in the \textit{Typical} scenario, cumulative profits for both policies tend to even out by the end of the optimization horizon.

%As profits are heavily dependent on the sample path, Table \ref{table:Model Table} reports the mean percentage of the posterior optimal solution under each modeling assumption. The standard deviation and worst case scenario are also displayed. Notice that the mean minus two standard deviations (an approximate 95 percent confidence interval if these percentages were assumed to be normally distributed), is highest in each case for the crossing state HSMM model. In addition, the worst performance relative to the posterior optimal is higher for the HSMM than those realized by the other two policies. These are good indicators that HSMM model produces policies that are not only effective, but robust as wel

\section{Conclusions}
\label{Conclusion}

In this paper we propose the use of crossing state hidden semi-Markov models for the stochastics involved in a common energy storage problem. The crossing state HSMM is shown to produce sample paths that replicate the crossing times of wind power forecast errors more accurately than common time series models. We show empirically that properly modeling this behavior leads to more robust control policies than when an AR(1) model, which consistently underestimates crossing time duration, is assumed. Though the HSMM-based policy is more robust in worst case scenarios, it is also as effective as the AR(1)-based policy in typical scenarios. A system operator could thus reduce risk without sacrificing expected profit simply by using this higher quality stochastic model.

The crossing state HSMM does introduce some complexity when designing a VFA-based policy. The partial unobservability of the crossing state requires a Bayesian update at each time step, and, furthermore, computing value functions for each possible system state is computationally expensive. An algorithm that computes solutions efficiently without much reduction in solution quality is thus necessary to use HSMMs in this problem setting. Backward approximate dynamic programming is shown to produce solutions that achieve 95 percent of optimality or above in a variety of test cases, using a fraction of the time needed to compute the exact solutions (and with reduced memory storage requirements as well). Backward ADP also outperforms other common control policies including approximate policy iteration, a tuned policy function approximation, and deterministic lookahead policies.

Comparing the two backward ADP algorithms used in this paper, Algorithm \ref{alg:Sampled MDP}, which fits value functions to post-decision states in lookup table form, proves to be not only more effective, but faster in this problem setting due to the compact post-decision state space. Meanwhile, Algorithm \ref{alg:Lin VFA}, which uses parametric VFAs, provides consistent performance at different sampling rates and requires minimal memory storage space. Finally, we note that the the approximate policies computed with a backward ADP algorithm assuming the crossing state HSMM for wind power forecast errors are more effective than the exact solutions assuming the AR(1) model when tested on historical sample paths (see Section \ref{Model Selection}). This justifies the use of the more complex stochastic model, even if an approximate dynamic programming algorithm is necessary to efficiently compute the policy as a result.

\begin{APPENDIX}{Extending the Crossing State HSMM to Model Electricity Prices}

\label{Electricity Prices}

\begin{figure}[b]
\centering
\includegraphics[width=5.5 in, height=2.2 in]{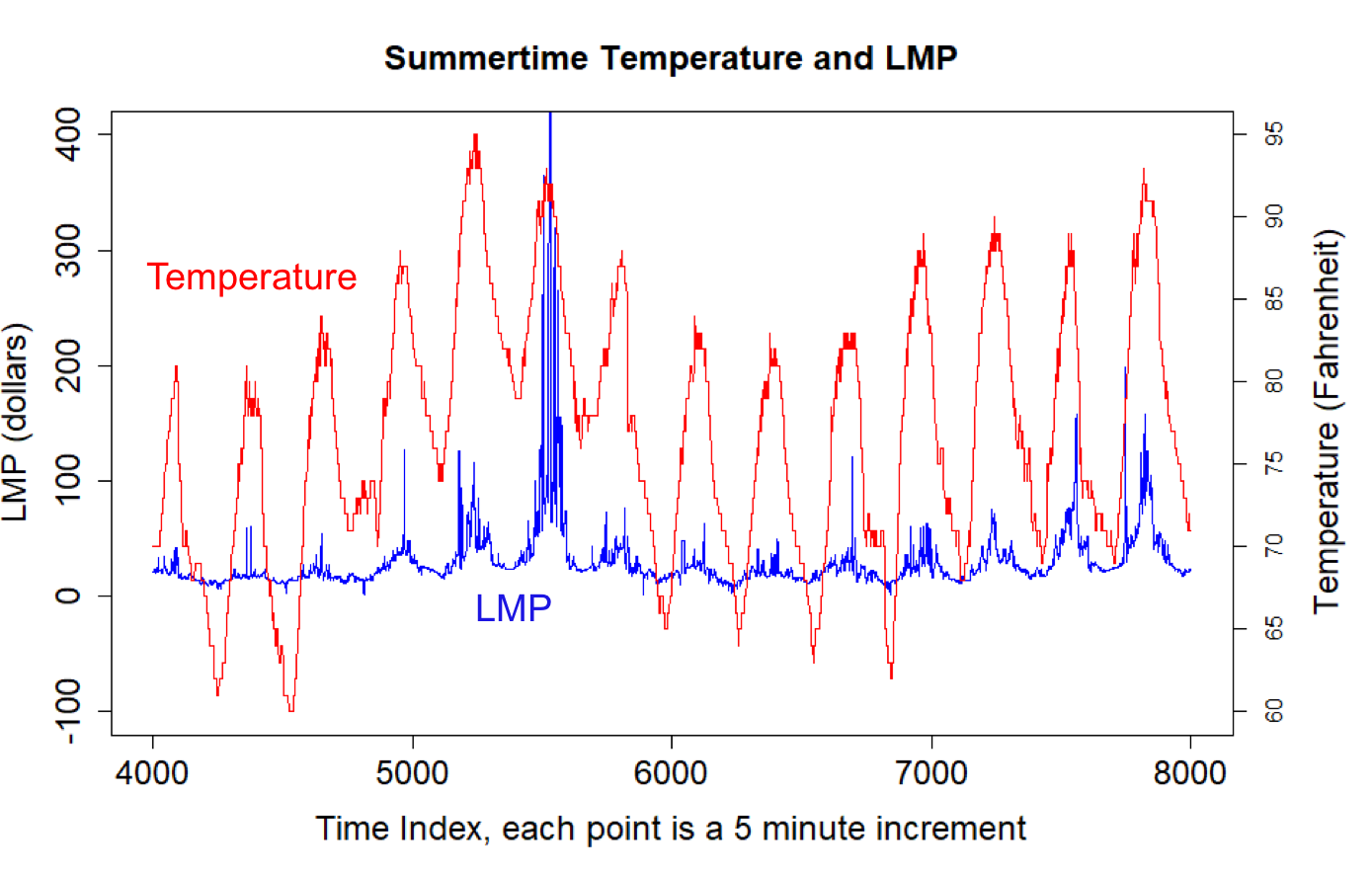}
\caption{LMP path and temperature data for the Princeton, NJ, area during two weeks in July 2015. LMPs spike with daily peaks and exhibit larger spikes when the peaks are higher.}
\label{LMPTempAct}
\end{figure}

As mentioned in Section \ref{WPM}, electricity prices exhibit much different behavior than wind power forecast errors, such as the tendency to spike under certain conditions, a heavy dependence on temperature, and a daily periodic pattern. These characteristics can be seen in Figure \ref{LMPTempAct}, in which a two week portion of LMP data from the Princeton, New Jersey area during the summer of 2015 is shown. Also shown is the observed temperature in the area during this period. We extend the crossing state model described in Sections \ref{WPM} and \ref{CompactInfo} to model summertime electricity prices by using a periodic reference series and conditioning price distributions on temperature.

We do not have access to external forecasts for LMPs here. Instead, we create our own forecast by letting $f_t^P=\mathbb{E}[P_t] \text{ } \forall t$ which is estimated by using the average price at each time of day observed in training data. This contrived forecast effectively captures the seasonal (daily) component of LMPs. As this reference series is deterministic over the optimization horizon, $\left\lbrace f^P_t \right\rbrace_{t=0}^T \in S_0$.

Observe from Figure \ref{LMPTempAct} that LMPs often spike during the high points of temperature each day. Additionally, the higher the peak temperature, the higher the spikes and price level in general. To capture this behavior, we incorporate temperature as an explanatory variable in the error generation portion of the model by conditioning the price distribution at time $t+1$ on the temperature forecast for time $t+1$. Ideally, this would be a forecast made at time $t$, but we assume a fixed temperature forecast over the optimization horizon as it is accurate enough to capture the general behavior of temperature in the near future and simplifies the optimization algorithm.

We first isolate two components of the temperature series -- the seasonal component and the trend component (identified by applying a length 50 moving average filter to the temperature series). Then, the model is conditioned on both of these variables. As both series take on continuous values, they must be aggregated into bins first. In this application, we use two bins for the seasonal component and three bins for the trend component. This granularity, along with the bin division points, were decided upon through trial and error.

Let $h^{s}_t$ be the seasonal component of the temperature series at time $t$. Let $h^{s,max}=\operatornamewithlimits{max}\limits_{t} h^{s}_t$. The explanatory variable $y^s_t$ may take on two values,
\begin{align}
y^s_t=\begin{cases}
&2 \text{ if } h^s_t \geq 0.75 \times h^{s,max}\\
&1 \text{ if } h^s_t < 0.75 \times h^{s,max}.
\end{cases}
\end{align}This division is intended to separate periods when the temperature is peaking from all other times of the day.

Let $h^{tr}_t$ be the trend component of the temperature series at time $t$. Let $h^{tr,max}=\operatornamewithlimits{max}\limits_{t}h^{tr}_t$. The explanatory variable $y^{tr}_t$ may take on three values,
\begin{align}
y^{tr}_t=\begin{cases}
&3 \text{ if } h^{tr}_t \geq 0.8 \times h^{tr,max}\\
&2 \text{ if } 0.3 \times h^{tr, max} \leq h^{tr}_t < 0.8 \times h^{tr,max}\\
&1 \text{ if } h^{tr}_t < 0.3 \times h^{tr,max}.
\end{cases}
\end{align} This division is intended to separate cool, average, and hot days.

Several LMP sample paths for a one day horizon are shown in Figure \ref{LMPSim} where the effect of conditioning on the variables $y_t^{s}$ and $y_t^{tr}$ can be observed. As $\left\lbrace y^s_t \right\rbrace_{t=0}^T$ and $\left\lbrace y^{tr}_t \right\rbrace_{t=0}^T$ are both deterministic series (stemming from the fact that we assume a fixed temperature forecast over the optimization horizon), these belong in the initial state $S_0$. 

\begin{figure}
\centering
\includegraphics[width=5.5 in, height=2.2 in]{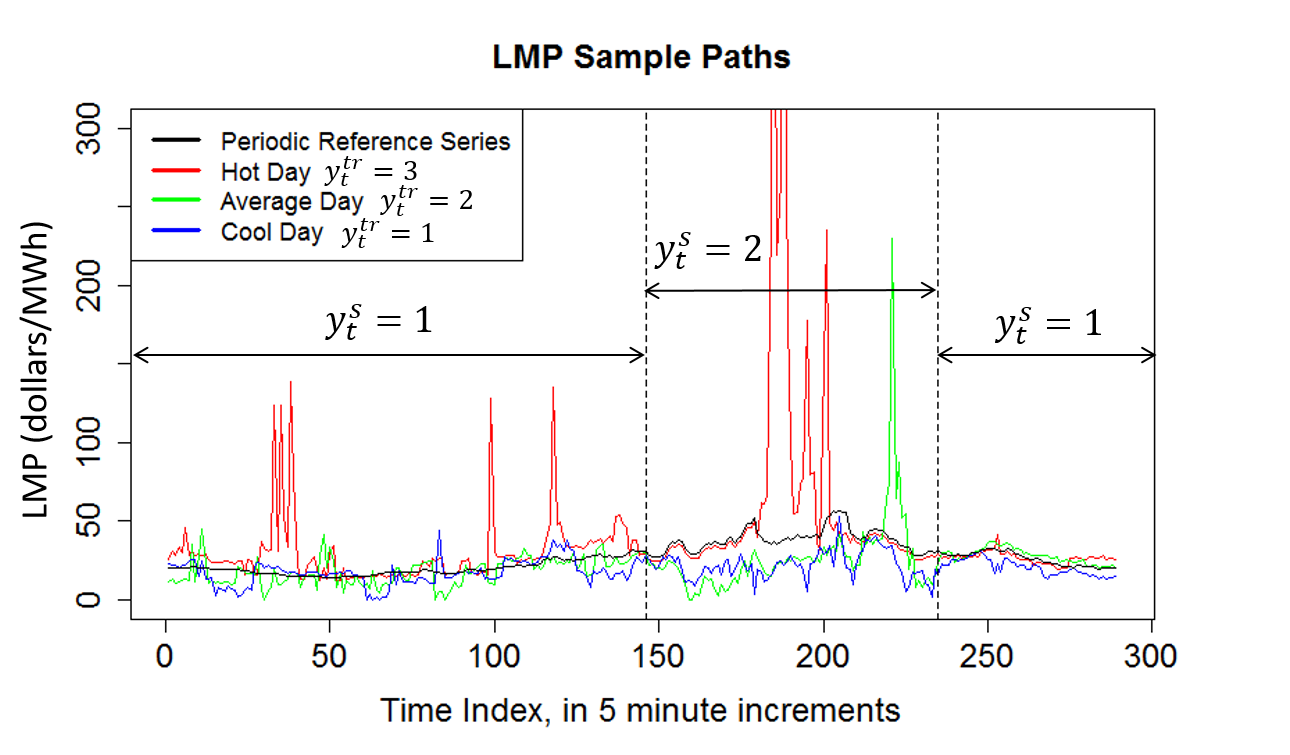}
\caption{Example LMP sample paths over different days. The effect of conditioning on temperature trend and seasonality separately can be observed.}
\label{LMPSim}
\end{figure}

With this additional conditioning, our error-from-reference density at time $t+1$ gives $\mathbf{P}(\hat{P}_{t+1}|I_t^{P}, y^s_{t+1}, y^{tr}_{t+1})$ for each information state $I_t^P$. Similar to the wind model, we also have a compact information state, $\tilde{I}_t^P=(C_t^P,D_t^P,\hat{P}^b_t)$, utilized when fitting value functions. The set of all compact information states, $\tilde{\mathcal{I}}^P$, has cardinality $2 \times m^P \times n^P$ as in the wind model.

As in Section \ref{Model Selection}, we compare the performance of the policies developed using this HSMM to model LMPs to those resulting from more common electricity price model assumptions. A very similar training and testing procedure is used here as in Section \ref{Model Selection}. Three models are used to train value functions (of lookup table form for post-decision states), of which the first two are commonly used to model electricity prices in the literature:
\begin{henumerate}
\item{A first-order Markov chain: $\mathbf{P}(\hat{P}_t|\hat{P}^b_{t-1})$ where the price errors are binned into discrete states to form conditional distributions of $\hat{P}_t$ given the bin of $\hat{P}_{t-1}$.}
\item{A model with mean reversion and jump-diffusion (MRJD): $\hat{P}_t=\gamma^P \hat{P}_{t-1}+\mathcal{N}(0,(s^{P})^2)+\mathbf{1}_{\left\lbrace U < p \right\rbrace}\mathcal{N}(0,(s^{J})^2)$, where $U \sim Unif[0,1]$, $s^{J}$ and $p$ are chosen to capture the standard deviation and frequency of the price spike events, and $s^{P}$ is the standard deviation of the data without the price spikes.}
\item{The crossing state model for electricity prices described in this paper.}
\end{henumerate} Note that we are modeling $\hat{P}_t$, the deviation from the time-dependent mean, not $P_t$ itself, such that each model has the same deterministic seasonality component. Additionally, the post-decision state space for each model is discretized such that each model has the same number of possible post-decision states. Finally, the same crossing state HSMM is used for the wind power forecast error process in all cases.

\begin{table}[]
\TABLE
{The effect of electricity price model choice on policy performance. Similar testing is performed as in Table \ref{table:Model Table}, except a slightly different set of \textit{Typical} and \textit{Worst Case} sample paths, described in the text, are used for policy evaluation.
\label{table:PModel Table}}
 % used for centering table
{\begin{tabular}{c || c c | c c | c c } % centered columns (4 columns)
\hline\hline 
\multirow{2}{*}{Electricity Price Model} & \multicolumn{2}{c}{Test Case 1 (Avg Temp)} & \multicolumn{2}{c}{Test Case 2 (Cool Temp)} & \multicolumn{2}{c}{Test Case 3 (Hot Temp)} \\
\hhline{~------}
  & Typical & Worst Case & Typical & Worst Case & Typical & Worst Case \\

\hline % inserts single horizontal line
 Markov &  \bf{-757} & \bf{-1045} & -313 & -609 & -2833 & -3383\\
  MRJD  & -849 & -1121 & -365 & -645 & -2772 & -3332 \\
  HSMM &  \bf{-757} & -1050  &\bf{-306} & \bf{-601} & \bf{-2687} & \bf{-3256}

 % is used to refer this table in the text
\end{tabular}}
{}
\end{table}

The test cases in Table \ref{table:PModel Table} are chosen such that the price models are tested on cool, average, and hot days. Historical price sample paths observed under each temperature condition are used for policy evaluation in the corresponding test case. These are paired with the \textit{Typical} and \textit{Worst Case} wind power sample paths described in Section \ref{Model Selection} to form $\Omega^{Typ}$ and $\Omega^{WC}$ here. Table \ref{table:PModel Table} reports $F^{\pi,Typ}$ and $F^{\pi,WC}$ for the VFA-based policies trained on each price model in each test case. Observe that the crossing state model for prices produces policies that perform similarly to policies based on more common stochastic models on average and cool temperature days, but are more robust when price spikes are larger and more frequent on hot days (whether wind has \textit{Typical} or \textit{Worst Case} behavior). Conditioning price distributions on a temperature forecast is clearly a key factor in this, as the model more accurately captures when price spikes occur, at what frequency they occur, and at what magnitude under several temperature conditions. The resulting policy based on this model can then plan an appropriate amount of storage based on this information.

 \end{APPENDIX}

% or


\section*{Acknowledgment}

The research was supported by NSF grant CCF-1521675 and DARPA grant FA8750-17-2-0027.

%\textsc{\ifCLASSOPTIONcaptionsoff
%  \newpage
%\fi}




\bibliographystyle{informs2014}
\bibliography{bibpaper}

\end{document}